\documentclass[aap,MSNbibl,dvips]{arximspdf}
\usepackage{graphicx}
\doi{10.1214/13-AAP970} %kopijuoti is PTS
\volume{24}
\issue{6}
\pubyear{2014}
\firstpage{2207}
\lastpage{2245}

\makeatletter
\newcommand{\widebar}{\overline}

\newcommand{\eqref}[1]{(\ref{#1})}
\newcommand{\bmu}{\bolds{\mu}}
\newcommand{\sX}{\mathsf{X}}
\newcommand{\bbb}{\mathbf{b}}
\newcommand{\bW}{\mathbf{W}}
\newcommand{\hZ}{\widehat{Z}}
\newcommand{\cB}{\mathcal{B}}

\newcommand{\bbE}{\mathbb{E}}
\newcommand{\bZero}{\mathbf{0}}
\newcommand{\bOne}{\mathbf{1}}
\newcommand{\beps}{\varepsilon}

\newcommand{\bphi}{\bolds{\phi}}
\newcommand{\bA}{\mathbf{A}}
\newcommand{\cC}{\mathcal{C}}
\newcommand{\bB}{\mathbf{B}}
\newcommand{\bD}{\mathbf{D}}
\newcommand{\hull}[1]{\langle{#1}\rangle}
\newcommand{\cK}{\mathcal{K}}
\newcommand{\cM}{\mathcal{M}}
\newcommand{\N}{{N}}
\newcommand{\bQ}{\mathbf{Q}}
\newcommand{\sS}{\mathcal{S}}
\newcommand{\T}{\fontsize{8.36pt}{9pt}\selectfont{{\textsf{T}}}}
\newcommand{\bu}{\mathbf{u}}
\newcommand{\bv}{\mathbf{v}}
\newcommand{\bL}{\mathbf{L}}
\newcommand{\bY}{\mathbf{Y}}
\newcommand{\bZ}{\mathbf{Z}}
\newcommand{\PRIMAL}{\mathrm{PRIMAL}}
\newcommand{\bSN}{\textbf{SN}}
\newcommand{\bBN}{\textbf{BN}}

\newcommand{\Prob}{\mathbb{P}}
\newcommand{\Reals}{\mathbb{R}}
\newcommand{\Integers}{\mathbb{Z}}
\newcommand{\IntegersP}{\mathbb{Z}_+}
\newcommand{\Naturals}{\mathbb{N}}
\newcommand{\RealsP}{\Reals_+}
\newcommand{\mmin}{\wedge}
\newcommand{\mmax}{\vee}
\newcommand{\bdot}{\cdot}

\newtheorem{proposition}[theorem]{Proposition}
\newtheorem{theorem}{Theorem}[section]
\newtheorem{corollary}[theorem]{Corollary}
\newtheorem{lemma}[theorem]{Lemma}
\newproclaim{definition}[theorem]{Definition}
\newproclaim{assumption}[theorem]{Assumption}

\def\green#1{{\color{green}#1}}
\def\old#1{}

\newcommand{\cI}{\mathcal{I}}
\newcommand{\bC}{\mathbf{C}}
\newcommand{\bM}{\mathbf{M}}
\newcommand{\bX}{\mathbf{X}}
\newcommand{\bR}{\mathbf{R}}
\newcommand{\be}{\mathbf{e}}

\newcommand{\bbI}{\mathbb{I}}
\newcommand{\bbZ}{\mathbb{Z}}
\makeatother

\begin{document}
\begin{frontmatter}

\title{Optimal queue-size scaling in switched networks}
\runtitle{Optimal scheduling}

\begin{aug}
\author[a]{\fnms{D.} \snm{Shah}\thanksref{t1}\ead[label=e1]{devavrat@mit.edu}},
\author[b]{\fnms{N. S.} \snm{Walton}\ead[label=e2]{n.s.walton@uva.nl}}
\and
\author[c]{\fnms{Y.} \snm{Zhong}\corref{}\ead[label=e3]{yz2561@columbia.edu}\thanksref{t1}}
\runauthor{D. Shah, N. S. Walton and Y. Zhong}
\affiliation{Massachusetts Institute of Technology, University of Amsterdam\\ and Columbia University}
\address[a]{D. Shah\\
Department of EECS\\
Massachusetts Institute of Technology\\
Cambridge, Massachusetts 02139\\
USA\\
\printead{e1}}
\address[b]{N. Walton\\
Korteweg-de Vries Institute for Mathematics\\
University of Amsterdam\\
1090 GE Amsterdam\\
The Netherlands\\
\printead{e2}}
\address[c]{Y. Zhong\\
Department of IEOR\\
Columbia University\\
New York, New York 10027\\
USA\\
\printead{e3}}%adresu isvedimo komanda gale!
\end{aug}
\thankstext{t1}{Supported by NSF TF collaborative project and NSF CNS CAREER project.
When this work was performed, the third author was affiliated with the
Laboratory for Information and Decision Systems as well as the
Operations Research Center at MIT. The third author is now affiliated
with the Department of Industrial Engineering and Operations Research
at Columbia University.}

%CCF-0728554, NSF CAREER project
%and XXX. Author's email addresses are {\tt devavrat@mit.edu}, {\tt
%jnt@mit.edu}, {\tt n.s.walton@uva.nl}, and

% HISTORY:
\received{\smonth{10} \syear{2011}}
\revised{\smonth{12} \syear{2012}}

% ABSTRACT
%
\begin{abstract}
We consider a switched (queuing) network in which there are constraints
on which queues may be served simultaneously; such networks have been
used to effectively model input-queued switches and wireless networks.
The scheduling policy for such a network specifies which queues to
serve at any point
in time, based on the current state or past history of the system.
In the main result of this paper, we provide a new class of
online scheduling policies that achieve optimal queue-size scaling
for a class of switched networks including input-queued switches.
In particular, it establishes the validity of a conjecture (documented
in Shah, Tsitsiklis and Zhong
[\textit{Queueing Syst.} \textbf{68} (2011) 375--384])
about optimal queue-size scaling for input-queued switches.
\end{abstract}

% KEYWORDS
\begin{keyword}[class=AMS]
\kwd{60K25}
\kwd{60K30}
\kwd{90B36}
\end{keyword}
\begin{keyword}
\kwd{Switched network}
\kwd{maximum weight scheduling}
\kwd{fluid models}
\kwd{state space collapse}
\kwd{heavy traffic}
\kwd{diffusion approximation}
\end{keyword}

\end{frontmatter}

\setcounter{footnote}{1}
%s1 #&#
\section{Introduction}\label{secintro}

A switched network consists of a collection of, say, $N$~queues,
operating in discrete time.
At each time slot, queues are offered service according to a \emph{service schedule} chosen
from a specified finite set, denoted by~$\sS$. The rule for choosing a
schedule from $\sS$ at
each time slot is called the \emph{scheduling policy}. New work may
arrive to each queue
at each time slot exogenously and work served from a queue may join
another queue or leave the network. We shall restrict our attention,
however, to the case
where work arrives in the form of unit-sized packets, and once it is
served from a
queue, it leaves the network, that is, the network is single-hop.

%
%Each queue has a dedicated exogenous arrival process, with specified
%mean
%arrival rates. Let $\lambda_i$ be the rate at which packets arrive to
%queue $i$ and
%let $\blambda= [\lambda_i]$. If schedule $\bolds{\sigma}\in\sS$, which
%also belongs to $\{0,1\}^N$,
%is chosen for servicing queues in a given time slot, and if $\sigma_i
%= 1$ then a
%packet leaves from queue $i$ in that time slot if queue $i$ were
%nonempty. We
%shall restrict to single-hop network, that is once a packet is served
%from a queue
%it leaves the network.

Switched networks are special cases of what Harrison \cite{harrisoncanonical,harrisoncanonicalcorr} calls
``stochastic processing networks.'' Switched networks are general
enough to model a variety of interesting applications. For example,
they have been used
to effectively model input-queued switches, the devices at the heart of high-end
Internet routers, whose underlying silicon architecture imposes
constraints on
which traffic streams can be transmitted simultaneously \cite{daibala}.
They have also been used to model multihop wireless networks in which
interference limits the amount of service that can be given to each
host \cite{tassiula1}. Finally, they can be instrumental in finding
the right
operational point in a data center~\cite{SWo}.

In this paper, we consider \textit{online} scheduling policies,
that is, policies that only utilize historical information (i.e., past
arrivals and scheduling decisions). The performance objective of interest
is the total queue size or total number of packets waiting to be served
in the network on average (appropriately defined). The questions that we
wish to answer are: (a) what is the minimal value of the performance objective
among the class of online scheduling policies,
and (b) how does it depend on the network structure, $\sS$,
as well as the effective load.

Consider a work-conserving $M/D/1$ queue with a unit-rate server in
which unit-sized
packets arrive as a Poisson process with rate $\rho\in(0, 1)$. Then,
the long-run average
queue-size scales\footnote{In this paper, by scaling of quantity we
mean its dependence
(ignoring universal constants) on~$\frac{1}{1-\rho}$ and/or the
number of queues, $N$,
as these quantities become large. Of particular interest is the scaling
of $\rho\to1$ and
$N\to\infty$, in that order. } as $1/(1-\rho)$. Such scaling
dependence of the average
queue size on $1/(1-\rho)$ (or the inverse of the \textit{gap}, $1-\rho
$, from the load to the
capacity) is a universally observed behavior in a large class of
queuing networks.
In a switched network, the scaling of the average total queue size
ought to
depend on the number of queues, $N$. For example,
consider $N$ parallel $M/D/1$ queues as described above.
Clearly, the average total queue size will scale as \mbox{$N/(1-\rho)$}. On
the other
hand, consider a variation where all of these queues pool their
resources into a single
server that works $N$ times faster. Equivalently, by a time change,
let each of the $N$ queues receive
packets as an independent Poisson process of rate $\rho/N$, and each
time a
common unit-rate server serves a packet from one of the nonempty
queues. Then, the
average total queue-size scales as $1/(1-\rho)$. Indeed, these are
instances of switched networks that differ in their scheduling set $\sS
$, which leads to
different queue-size scalings. Therefore, a natural question is the
determination of
queue-size scaling in terms of $\sS$ and $(1-\rho)$, where $\rho$ is
the effective load. In the context of an $n$-port input-queued switch with
$N = n^2$ queues, the optimal scaling of average total queue size has been
conjectured to be $n/(1-\rho)$, that is, $\sqrt{N}/(1-\rho)$ \cite{STZopen}.

As the main result of this paper, we propose a new online scheduling
policy for
any single-hop switched network. This policy effectively emulates an
insensitive bandwidth
sharing network with a product-form stationary distribution with each
component of this
product-form behaving like an $M/M/1$ queue. This crisp description of
stationary distribution allows us to
obtain precise bounds on the average queue sizes under this policy.
This leads to establishing, as a corollary of our result, the validity
of a conjecture stated in \cite{STZopen} for input-queued switches.
In general, it provides explicit bounds on the average total queue size
for any switched network. Furthermore, due to the explicit bound on the
stationary
distribution of queue sizes under our policy, we are able to establish
a form
of large-deviations optimality of the policy for a large class of
single-hop switched networks,
including the input-queued switches, and the independent-set model of
wireless networks,
when the underlying interference graph is bipartite, for example, and
more generally, perfect.

The conjecture from \cite{STZopen} that we settle in this paper,
states that
in the heavy-traffic regime (i.e., $\rho\rightarrow1$), the optimal
average total queue-size scales as \mbox{$\sqrt{N}/(1-\rho)$}.
%We note that the validity of the conjecture in \cite{STZopen} for
%input-queued
%switches, stating that optimal average total queue-size scales as $
The validity of this conjecture
is a significant improvement over the best-known bounds of $O(N/(1-\rho))$
(due to the moment bounds of \cite{MT93} for the maximum weight
policy) or
$O(\sqrt{N}\log N/(1-\rho)^2)$ (obtained by using a batching policy
\cite{NeelyModiano}).

% when $1/(1-\rho)$ scales faster than $n^2$.
Our analysis consists of two principal components. First, we
propose and analyze a scheduling mechanism that is able to emulate, in
discrete time, \textit{any} continuous-time bandwidth allocation within a
bounded degree of error. This scheduler maintains a continuous-time
queuing process and tracks its own queue size process. If, valued under
a certain decomposition, the gap between the idealized continuous-time
process and the real queuing process becomes too large, then an
appropriate schedule is allocated. Second, we implement specific
bandwidth allocation named the store-and-forward allocation policy
(SFA). This policy was first considered by Massouli\'e, and was consequently
discussed in the thesis of Prouti\`ere \cite{PTh}, Section~3.4. It was
shown to
be insensitive with respect to phase-type service distributions in
works by
Bonald and Prouti\`ere \cite{bonaldproutiere1,bonaldproutiere2}. The
insensitivity of this
policy for general service distributions was established by Zachary~\cite{zachary}. The store-and-forward policy is closely related to
the classical product-form multi-class queuing network, which have
highly desirable queue-size scalings. By emulating these queuing
networks, we are able to translate results which render optimal
queue-size bounds for a switched network. An interested reader is
referred to \cite{walton} and \cite{KMW} for an in-depth discussion
on the relation between this policy, the proportionally fair
allocation, and multi-class queuing networks.

%s1.1 #&#
\subsection{Organization} In Section~\ref{secmodel}, we specify a
stochastic switched network model. In Section~\ref{secrelated}, we
discuss related works. Section~\ref{secinsensitive} details the
necessary background on the insensitive store-and-forward
bandwidth allocation (SFA) policy. The main result of the paper
is presented and proved in Section~\ref{secmain}.
We first describe the policy for single-hop switched networks,
and state our main result, Theorem~\ref{thmmain}.
This is followed by a discussion of the optimality of the policy.
We then provide a proof of Theorem~\ref{thmmain}. %First, we describe
%the
%policy and its performance in the context of input-queued switches.
%This is followed by the policy and performance
%analysis for single-hop switched networks.
A discussion of
directions for future work is provided in Section~\ref{secconc}.

\subsubsection*{Notation}
Let $\Naturals$ be the set of natural numbers $\{1,2,\ldots\}$,
let $\IntegersP=\{0,1,\break 2,\ldots\}$, let $\Reals$ be the set of real numbers
and let $\RealsP=\{x\in\Reals\dvtx x\geq0\}$.
Let $\bbI[A]$ be the indicator function of an event $A$, %where
%$1_{true}=1$ and $1_{\mathrm{false}}=0$.
Let $x\mmin y=\min(x,y)$, $x\mmax y=\max(x,y)$ and $[x]^+=x\mmax0$.
When $\mathbf{x}$
is a vector, the maximum is taken componentwise.

We will reserve bold letters for vectors in
$\Reals^\N$, where $N$ is the number of queues. For example,
$\mathbf{x}=[x_n]_{1\leq n\leq N}$.
Superscripts on vectors are used to denote labels, not exponents,
except where otherwise noted;
thus, for example, $(\mathbf{x}^0,\mathbf{x}^1,\mathbf{x}^2)$ refers
to three arbitrary
vectors.
Let $\bZero$ be the vector of all 0s and $\bOne$ the
vector of all~1s. The vector $\be_i$ is the $i$th
unit vector, with all components being $0$ but the $i$th component
equal to $1$.
We use the norm $|\mathbf{x}|=\max_n|x_n|$. For vectors $\bu$ and
$\bv$, we let
$\bu\bdot\bv= \sum_{{n}=1}^\N u_{n}v_{n}$.
%and functions $f\dvtx \Reals\to\Reals$, we let
%f(\bu) = \bigl[ f(u_\n) \bigr]_{1\leq\n\leq\N},
%and let matrix multiplication take precedence over dot product so that
Let $\bA^{\T}$ be the transpose of matrix $\bA$.
For a set $\sS\subset\Reals^N$, denote its convex hull by $\hull
{\sS}$.
For $n \in\Naturals$, let $n! = \prod_{\ell= 1}^n \ell$ be the
factorial of $n$,
and by convention, $0! = 1$.

%s2 #&#
\section{Switched network model}\label{secmodel}
We now introduce the switched network model. Section~\ref{secmodelqueue}
describes the general system model, Section~\ref{secmodelstochastic}
lists the
probabilistic assumptions about the arrival process and Section~\ref{secmodelquantity}
introduces some useful definitions.

%s2.1 #&#
\subsection{Queueing dynamics}\label{secmodelqueue}
Consider a collection of $N$ queues. Let time be discrete, and indexed by
$\tau\in\{0,1,\ldots\}$. Let $Q_i(\tau)$ be the amount of work in
queue $i\in\{1,\ldots,N\}$ at time slot $\tau$.
Following our general notation for vectors, we write
$\bQ(\tau)$ for $[Q_i(\tau)]_{1\leq i\leq N}$. The initial queue
sizes are $\bQ(0)$. Let $A_i(\tau)$ be the total amount of
work arriving to queue $i$, and $B_i(\tau)$ be the cumulative
potential service to queue $n$, up to time $\tau$, with
$\bA(0)=\bB(0)=\bZero$.

We first define the queuing dynamics for a single-hop switched network.
Defining $d\bA(\tau)=\bA(\tau+1)-\bA(\tau)$ and
$d\bB(\tau)=\bB(\tau+1)-\bB(\tau)$, the basic Lindley recursion
that we will consider is
%
%e1 #&#
\begin{equation}
\label{eqlindley1} \bQ(\tau+1) = \bigl[ \bQ(\tau)-d\bB(\tau) \bigr]^+ +
d\bA(\tau),
\end{equation}
where the operation $[\cdot]^+$ is applied componentwise. The fundamental
switched network constraint is that there is some finite set
$\sS\subset\RealsP^N$ such that
%
%e2 #&#
\begin{equation}
\label{eqscheduling} d\bB(\tau)\in\sS\qquad\mbox{for all }\tau.
\end{equation}
For the purpose of this work, we shall focus on $\sS\subset\{0,1\}^\N$.
We will refer to $\bolds{\sigma}\in\sS$ as a schedule and $\sS$ as the
set of
allowed schedules. In the applications in this paper, the schedule is
chosen based on current queue sizes, which is why it is
natural to write the basic Lindley recursion as \eqref{eqlindley1}
rather than the more standard $[\bQ(\tau)+d\bA(\tau)-d\bB(\tau)]^+$.

For the analysis in this paper, it is useful to keep track of two other
quantities.
Let $Z_i(\tau)$ be the cumulative amount of idling at
queue $n$, defined by $\bZ(0)=\bZero$ and
%
%e3 #&#
\begin{equation}
\label{eqdefidling} d\bZ(\tau) = \bigl[d\bB(\tau)-\bQ(\tau) \bigr]^+,
\end{equation}
where $d\bZ(\tau)=\bZ(\tau+1)-\bZ(\tau)$. Then, \eqref
{eqlindley1} can be rewritten as
%
%e4 #&#
\begin{equation}
\label{eqdiscretequeuesinglehop} \bQ(\tau) = \bQ(0) + \bA(\tau) - \bB
(\tau) + \bZ(\tau).
\end{equation}
Also, let $S_{\bolds{\sigma}}(\tau)$ be the cumulative amount of time
that is spent on using schedule~$\bolds{\sigma}$ up to time $\tau$, so that
%
%e5 #&#
\begin{equation}
\label{eqservice} \bB(\tau)=\sum_{\bolds{\sigma}\in\sS} S_{\bolds{\sigma}}(
\tau) \bolds{\sigma}.
\end{equation}

A policy that decides which schedule to choose at each time slot
$\tau\in\IntegersP$ is called a \emph{scheduling policy}.
In this paper, we will be interested in online scheduling policies.
That is, the scheduling decision at time $\tau$ will be based on
historical information, that is, the cumulative arrival process $\bA
(\cdot)$
%and scheduling decisions $\bB(\cdot)$ till time $\tau$.
till time $\tau$.

%s2.2 #&#
\subsection{Stochastic model}\label{secmodelstochastic}

We shall assume that the exogenous arrival process for each
queue is independent and Poisson. Specifically, unit-sized packets
arrive to queue $i$ as a Poisson process of rate $\lambda_i$.
Let $\bolds{\lambda}= [\lambda_i]_{i=1}^N$ denote the vector of
all arrival rates.
The results presented in this paper extend to more general
arrival process with i.i.d. interarrival times with finite means,
using a \textit{Poissonization} trick.
We discuss this extension in Section~\ref{secconc}.

%s2.3 #&#
\subsection{Useful quantities}\label{secmodelquantity}

We shall assume that the scheduling constraint set $\sS$ is \textit{monotone}.
This is captured in the following assumption.
%
%as2.1 #&#
\begin{assumption}[(Monotonicity)]\label{assmonotone}
If $\sS$ contains a schedule, then $\sS$ also contains all of
its sub-schedules. Formally, for any $\bolds{\sigma}\in\sS$,
if $\bolds{\sigma}' \in\{0,1\}^N$ and $\bolds{\sigma}'\leq\bolds{\sigma
}$ componentwise,
then $\bolds{\sigma}' \in\sS$.
\end{assumption}
Without loss of generality, we will assume that each unit vector $\be_i$
belongs to $\sS$.
Next, we define some quantities that will be useful in the remainder
of the paper.
%
%de2.2 #&#
\begin{definition}[(Admissible region)]\label{dfadmissible}
Let $\sS\subset\{0,1\}^N$ be the set of allowed schedules. Let
$\hull{\sS}$ be the convex hull of $\sS$, that is,
\begin{eqnarray*}
\hull{\sS} &=& \biggl\{ \sum_{\bolds{\sigma}\in\sS}
\alpha_{\bolds{\sigma}} \bolds{\sigma}\dvtx  \sum_{\bolds{\sigma}\in\sS}
\alpha_{\bolds{\sigma}} =1 \mbox{ and }\alpha_{\bolds{\sigma}}
\geq0, \mbox{ for all }\bolds{\sigma} \biggr\}.
\end{eqnarray*}
Define the \emph{admissible region} $\cC$ to be
\begin{eqnarray*}
\cC&=& \bigl\{\bolds{\lambda}\in\RealsP^N \dvtx  \bolds{\lambda}\leq
\bolds{\sigma}\mbox{ componentwise, for some }\bolds{\sigma}\in
\hull{\sS}
\bigr\}.
\end{eqnarray*}
\end{definition}
Note that under Assumption~\ref{assmonotone},
the capacity region $\cC$ and the convex hull $\hull{\sS}$
of $\sS$ coincide.

Given that $\hull{\sS}$ is a polytope contained
in $[0,1]^N$, there exists an integer $J \geq1$, a matrix
$\bR\in\RealsP^{J \times N}$ and a vector $\bC\in\RealsP^{J}$
such that
%
%e6 #&#
\begin{eqnarray}
\label{eqrankeq} \hull{\sS} & =& \bigl\{ \mathbf{x}\in[0,1]^N \dvtx  \bR
\mathbf{x}\leq\bC\bigr\}.
\end{eqnarray}
We call $J$ the \textit{rank} of $\hull{\sS}$ in the representation
\eqref{eqrankeq}.
When it is clear from the context, we simply call $J$ the rank of
$\hull{\sS}$.
Note that this rank may be different from the rank of matrix $\bR$.
Our results will exploit the fact that the rank $J$ may be an order of
magnitude smaller than $N$.

%
%de2.3 #&#
\begin{definition}[(Static planning problems and load)]\label{dfstaticplanning}
Define the static planning optimization problem $\PRIMAL(\bolds
{\lambda})$ for
$\bolds{\lambda}\in\RealsP^N$ to be
%
%e7 #&#
%e8 #&#
%e9 #&#
\begin{eqnarray}
\label{defprimal} \mbox{minimize}&\qquad& \sum_{\bolds{\sigma}\in\sS}
\alpha_{\bolds{\sigma}},
\\
\mbox{subject to}&\qquad& \bolds{\lambda}\leq\sum_{\bolds{\sigma}\in\sS}
\alpha_{\bolds{\sigma}} \bolds{\sigma}, \label{constreqprimal}
\\
&\qquad&\alpha_{\bolds{\sigma}}\in\RealsP\qquad\mbox{for all }\bolds
{\sigma} \in
\sS.
\end{eqnarray}
Define the induced load by $\bolds{\lambda}$, denoted by $\rho
(\bolds{\lambda})$,
as the value of the optimization problem $\PRIMAL(\bolds{\lambda})$.
%For any admissible $\blambda$, $\rho(\blambda) \leq1$.
\end{definition}
Note that $\bolds{\lambda}$ is admissible if and only if $\rho
(\bolds{\lambda}) \leq1$.
It also follows immediately from Definition~\ref{dfstaticplanning} that
%
%e10 #&#
\begin{equation}
\label{eqdef-load-2} \rho(\bolds{\lambda}) = \inf\{\gamma\geq0\dvtx  \bR
\bolds{\lambda}
\leq\gamma\bC\},
\end{equation}
and $\bolds{\lambda}$ is admissible if and only if $\bR
\bolds{\lambda}\leq\bC$, componentwise.

In the sequel, we will often consider the quantities
$\tilde{\rho}_j = \sum_i R_{ji}\lambda_i / C_j$, for $j \in\{1, 2,\ldots, J\}$,
which can be interpreted as loads on individual ``resources'' of the system
(this interpretation will be made precise in Section~\ref{secinsensitive}).
They are closely related to the system load $\rho(\bolds{\lambda})$.
We formalize this relation in the following lemma,
whose proof is straightforward and omitted.
%
%le2.4 #&#
\begin{lemma}\label{lemsys-ind-load}
Consider a nonnegative matrix $\bR\in\Reals_+^{J\times N}$
and a vector $\bC\in\Reals^J$ with $C_j > 0$ for all $j$.
For a nonnegative vector $\bolds{\lambda}\in\RealsP^N$,
define $\rho(\bolds{\lambda})$ by \eqref{eqdef-load-2}
and $\tilde{\rho}_j = (\sum_i R_{ji}\lambda_i ) / C_j$.
Then $\rho(\bolds{\lambda}) = \max_j \tilde{\rho}_j$.
\end{lemma}

The following is a simple and useful property of $\rho(\cdot)$: for
any $\mathbf{a}, \bbb\in\RealsP^N$,
%
%e11 #&#
\begin{eqnarray}
\label{eqrhosubadd} \rho(\mathbf{a}+ \bbb) & \leq&\rho(\mathbf{a}) + \rho
(\bbb).
\end{eqnarray}

%%%\iffalse
%%%%
%%%\begin{assumption}\label{assumption}
%%%By definition, the admissible region $\cC$ is a polytope contained
%%%in $[0,1]^N$. Let represented as $\{ \mathbf{x}\in[0,1]^N \dvtx  \bR
%%%\mathbf{x}\leq\bC\}$
%%%where $\bR\in\RealsP^{J \times N}$ matrix and $\bC\in\RealsP^J$
%%%for some $J \leq N$.
%%%We shall assume that for any $\mathbf{x}\in\cC$, there exists $\hat
%%%{\mathbf{x}}$ so that
%%%$\bR\hat{\mathbf{x}} = \bC$ and $\hat{\mathbf{x}} \geq\mathbf
%%%{x}$ componentwise.
%%%\end{assumption}
%%%%
%%%\fi

%s2.4 #&#
\subsection{Motivating example}\label{sseciqintro}

An Internet router has several input ports and output ports. A data
transmission cable is attached to each of these ports. Packets arrive
at the input ports. The function of the router is to work out which
output port each packet should go to, and to transfer packets to the
correct output ports. This last function is called \emph{switching}.
There are a number of possible switch architectures; we will consider
the commercially popular input-queued switch architecture.

%
%f1 #&#
\begin{figure}%[b]

\includegraphics{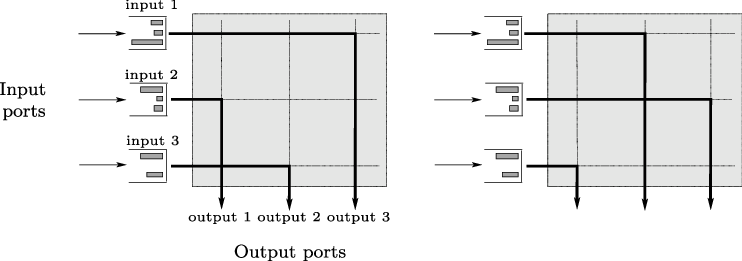}

\caption{An input-queued switch, and two example matchings of inputs
to outputs.}\label{figswitch-iq}
\end{figure}

Figure~\ref{figswitch-iq} illustrates an input-queued switch with
three input ports and three output ports. Packets arriving at input
$k$ destined for output $\ell$ are stored at input port~$k$, in queue
$Q_{k,\ell}$, thus there are $N=9$ queues in total. (For this example,
it is more natural to use double indexing, e.g., $Q_{3,2}$, whereas
for general switched networks it is more natural to use single indexing,
e.g., $Q_i$ for $1\leq i\leq N$.)

The switch operates in discrete time. At each time slot, the switch
fabric can transmit a number of packets from input ports to output
ports, subject to the two constraints that each input can transmit
at most one packet, and that each output can receive at most one
packet. In other words, at each time slot the switch can choose a
\emph{matching} from inputs to outputs. The schedule
$\bolds{\sigma}\in\RealsP^{3\times3}$ is given by $\sigma_{k,\ell}=1$ if
input port $k$ is matched to output port $\ell$ in a given time slot,
and $\sigma_{k,\ell}=0$ otherwise. The matching constraints require
that $\sum_{m=1}^3 \sigma_{k,m} \leq1$ for $k = 1, 2, 3$,
and $\sum_{m=1}^3 \sigma_{m, \ell} \leq1$ for $\ell= 1, 2, 3$.
Figure~\ref{figswitch-iq} shows two possible matchings. On the
left-hand side, the matching allows a packet to be transmitted
from input port 3 to output port 2, but since $Q_{3,2}$ is
empty, no packet is actually transmitted.

In general, for an $n$-port switch, there are $N = n^2$ queues.
The corresponding schedule set $\sS$ is defined as
%
%e12 #&#
\begin{equation}
\label{eqiq-sched-set} \sS= \Biggl\{\bolds{\sigma}\in\{0,1\}^{n\times
n} \dvtx  \sum
_{m=1}^n \sigma_{k,m} \leq1,
\sum_{m=1}^n \sigma_{m, \ell} \leq1,
1\leq k, \ell\leq n \Biggr\}.
\end{equation}
It can be checked that $\sS$ is \textit{monotone}. Furthermore, due to
Birkhoff--von Neumann theorem, \cite{BrK,VN}, the
convex hull of $\sS$ is given by
%
%e13 #&#
\begin{equation}
\label{eqiq-adm-region} \hull{\sS} = \Biggl\{\mathbf{x}\in
[0,1]^{n\times n} \dvtx  \sum
_{m=1}^n x_{k,m} \leq1, \sum
_{m=1}^n x_{m, \ell} \leq1, 1\leq k,
\ell\leq n \Biggr\}.
\end{equation}
Thus, the rank of $\hull{\sS}$ is less than or equal to $2n = 2\sqrt
{N}$ for an $n$-port switch.
Finally, given an arrival rate matrix\footnote{Not a vector, for
notational convenience, as discussed
earlier.} $\bolds{\lambda}\in[0,1]^{n\times n}$, $\rho
(\bolds{\lambda})$ is given by
\[
\rho(\bolds{\lambda}) = \max_{1\leq k, \ell\leq n} \Biggl\{\sum
_{m=1}^n \lambda_{k,m}, \sum
_{m=1}^n \lambda_{m,\ell} \Biggr\}.
\]

%s3 #&#
\section{Related works}\label{secrelated}

The question of determining the optimal scaling of queue sizes in switched
networks, or more generally, stochastic processing networks, has been an
important intellectual pursuit for more than a decade. The complexity
of the generic stochastic processing network makes this task extremely
challenging. Therefore, in search of tractable analysis, most of the
prior work has been on trying to understand optimal scaling and scheduling
policies for \textit{scaled} systems: primarily, with respect to fluid and
heavy-traffic scaling, {that is, $\rho\to1$.}

In heavy-traffic analysis, one studies the queue-size behavior under
a diffusion (or heavy-traffic) scaling. This regime was first considered
by Kingman \cite{kingmanht}; since then, a substantial body of theory has
developed, and modern treatments can be found in \cite
{mike2,bramson,williams,whittspl}. Stolyar \cite{stolyar} has studied
a class of myopic
scheduling policies, known as the maximum weight policy,
introduced by Tassiulas and Ephremides \cite{tassiula1}, for a generalized
switch model in the diffusion scaling. In a general
version of the maximum weight policy, a schedule with maximum weight
is chosen at each time step, with the weight of a schedule being equal
to the
sum of the weights of the queues chosen by that schedule. The weight of
a queue is a function of its size. In particular, for the choice of
one parameter
class of functions parameterized by $\alpha> 0$, $f(x)=x^\alpha$, the
resulting
class of policies are called the maximum weight policies with parameter
$\alpha> 0$,
and denoted as MW-$\alpha$.

In \cite{stolyar}, a complete characterization of the diffusion
approximation for the
queue-size process was obtained, under a condition known
as ``\emph{complete resource pooling},'' when the network is operating
under the MW-$\alpha$
policy, for any $\alpha> 0$.
%This condition effectively requires that
%there exists a scheduling policy which is able to balance the weights
%of all the
%heavily loaded queues.
Stolyar \cite{stolyar} showed the remarkable result
that the limiting queue-size vector lives in a one-dimensional state space.
Operationally, this means that all one needs to keep track of is the
one-dimensional total amount of work in the system (called the
\emph{rescaled workload}), and at any point in time one can assume that
the individual queues have all been balanced. Furthermore, it was established
that a max-weight policy minimizes the rescaled workload induced by any
policy under the
heavy-traffic scaling (with complete resource pooling).
Dai and Lin \cite{LD05,LD08} have established that a similar result holds
(with complete resource pooling) in the more general setting of
a stochastic processing network. In summary, under the complete
resource pooling condition, the results in \cite{stolyar,LD05,LD08} imply
that the performance of the maximum weight policy in an input-queued
switch, or more generally in a stochastic processing network, is always
optimal (in the diffusion limit, and when each queue size is
appropriately weighted).
These results suggest that the average total queue-size scales as
$1/(1-\rho)$
in the $\rho\to1$ limit. However, such analyses do not capture
the dependence on the network scheduling structure $\sS$. Essentially,
this is because the complete resource
pooling condition reduces the system to a one-dimensional space (which
may be highly dependent on a network's structure), and optimality
results are then initially expressed with respect to this
one-dimensional space.

Motivated to capture the dependence of the queue sizes on the network
scheduling structure $\sS$, a heavy-traffic analysis of switched
networks with
multiple bottlenecks (without resource pooling) was pursued by Shah and
Wischik \cite{SW}.
They established the so-called multiplicative state space collapse, and
identified a member, denoted by MW-$0^+$ (obtained by taking $\alpha
\to0$),
of the class of maximum-weight policies as optimal with respect to a critical
fluid model. In a more recent work, Shah and Wischik \cite{SWo}
established the optimality of
MW-$0^+$ with respect to overloaded fluid models as well. However, this
collection of works stops short of establishing optimality for diffusion
scaled queue-size processes.

Finally, we take note of the work by Meyn \cite{Meyn08}, which
establishes that
a class of generalized maximum weight policies achieve logarithmic [in
$1/(1-\rho)$]
regret with respect to an optimal policy under certain conditions.

In a related model---the bandwidth-sharing network model---Kang et~al.
\cite{kelly-williamsssc} have
established a diffusion approximation for the proportionally fair
bandwidth allocation policy,
assuming a technical ``local traffic'' condition, but without assuming
complete resource
pooling.\footnote{Kang et~al. \cite{kelly-williamsssc} assume that
critically loaded traffic is
such that all the constraints are saturated simultaneously.} They show
that the
resulting diffusion approximation has a product-form stationary distribution.
Shah, Tsitsiklis and Zhong \cite{STZ} have recently established that
this product-form stationary
distribution is indeed the limit of the stationary distributions of the original
stochastic model (an interchange-of-limits result). As a consequence,
if one could
utilize a scheduling policy in a switched network that corresponds to
the proportionally
fair policy, then the resulting diffusion approximation will have a product-form
stationary distribution, as long as the effective network scheduling structure
$\sS$ (precisely $\hull{\sS}$) satisfies the ``local traffic condition.''
Now, proportional fairness is a continuous-time rate allocation policy
that usually
requires rate allocations that are a convex combination of multiple schedules.
In a switched network, a policy must operate in discrete time and
has to choose one schedule at any given time from a finite discrete set
$\sS$.
For this reason, proportional fairness cannot be implemented directly.
However, a natural
randomized policy inspired by proportional fairness is likely to have the
same diffusion approximation (since the fluid models would be
identical, and the entire
machinery of Kang et~al. \cite{kelly-williamsssc}, building upon the
work of Bramson \cite{bramson} and Williams \cite{williams},
relies on a fluid model). As a consequence, if $\sS$ (more accurately,
$\hull{\sS}$)
satisfies the ``local traffic condition,'' then effectively the
diffusion-scaled queue sizes
would have a product-form stationary distribution, and would result in bounds
similar to those implied by our results. In comparison, our results
are nonasymptotic,
in the sense that they hold for any admissible load, have a
product-form structure, and do not require technical
assumptions such as the ``local traffic condition.'' Furthermore, such
generality is needed because there are popular
examples, such as the input-queued switch, that do \textit{not} satisfy the
``local traffic condition.''

Another line of works---so-called large-deviations analysis---concerns
exponentially decaying bounds on the tail probability
of the steady-state distributions of queue sizes.
Venkataramanan and Lin \cite{VL-LDP} established
that the maximum weight policy with weight parameter $\alpha> 0$,
MW-$\alpha$, optimizes the tail exponent of the $1+\alpha$ norm of the
queue-size vector.
Stolyar \cite{Stolyar-LDP} showed that a so-called ``exponential rule''
optimizes the tail exponent of the max norm of the queue-size vector.
However, these works do not characterize the tail
exponent explicitly. See \cite{STZSIGM} which has the best-known \textit{explicit} bounds on the tail exponent.

In the context of input-queued switches, the example that has
primarily motivated this work, the policy that we propose has the
average total queue size bounded within factor $2$ of the same quantity
induced by \emph{any} policy, in the heavy-traffic limit.
Furthermore, this result does not require conditions like complete
resource pooling.
More generally, our policy provides nonasymptotic bounds on queue sizes
for every arrival rate and switch size. The policy even admits
exponential tail bounds with respect to the stationary distribution,
and the exponent of these tail bounds is \emph{optimal}. These results
are significant improvements to the state-of-the-art bounds for best
performing policies for input-queued switches.
As noted in the \hyperref[secintro]{Introduction}, our bound on the average total queue size
is $\sqrt{N}$ times better than the existing bound for the
maximum-weight policy,
and $\log N/(1-\rho)$ times better than that for the batching policy
in \cite{NeelyModiano}.
(Here $N$ is the number of queues, and $\rho$ the system load.)
For further details of these results, see \cite{STZopen}.

For a generic switched network, our policy induces average
total queue size that scale linearly with the \textit{rank} of
$\hull{\sS}$, {under the diffusion scaling}. % and the usual $1/(1-
This is in contrast to the best-known bounds, such as those for
maximum weight policy, where the average queue-size scales as $N$,
under the diffusion scaling. %$N/(1-\rho)$.
Therefore, whenever the \textit{rank} of $\hull{\sS}$ is smaller than
$N$ (the
number of queues), our policy provides tighter bounds. Under our policy,
queue sizes admit exponential tail bounds. The bound on the distribution
of queue sizes under our policy leads to an explicit characterization of
the tail exponent, which is optimal for
{a wide range} of single-hop switched networks,
including input-queued switches and the independent-set model of
wireless networks,
{when the underlying interference graph is perfect}.

%any single-hop switched network.

%s4 #&#
\section{Insensitivity in stochastic networks}\label{secinsensitive}

This section recalls the background on insensitive stochastic networks that
underlies the main results of this work. We shall focus on
descriptions of the insensitive bandwidth allocation in so-called
bandwidth-sharing networks operating in continuous time. %in Section
Properties of these insensitive networks are provided in the \hyperref[app]{Appendix}.
%Justifications of claims made in this section are provided in the
%Appendix.

%Here we will see that bandwidth-sharing networks
%can be regarded as special instances of PS networks described in
%Section~\ref{ssecpsn}.
%In this context, rates are designated bandwidths,
%and bandwidth allocations are subject to capacity constraints.

We consider a bandwidth-sharing network operating in
continuous time with capacity constraints. The particular bandwidth-sharing
policy of interest is the \mbox{store-and-forward}
allocation (SFA) mentioned earlier.
We shall use the SFA as an idealized policy to design online scheduling
policies for switched networks. We now describe the precise model,
the SFA policy, and its performance properties.

\subsection*{Model} Let time be continuous and indexed
by $t \in\RealsP$. Consider a network with $J \geq1$ resources
indexed from $1,\ldots, J$. Let there be $N$ routes,
and suppose that each \emph{packet} on route $i$
consumes an amount $R_{ji}\geq0$ of resource $j$, for each $j \in\{1,
2, \ldots, J\}$. Let $\cK$ be the set of all resource--route pairs
$(j, i)$ such that route $i$ uses resource $j$, that is,
$\cK= \{(j,i) \dvtx  R_{ji} > 0\}$.
Without loss of generality, we assume that for each $i \in\{1, 2,\ldots, N\}$,
$\sum_{j=1}^J R_{ji} > 0$.
Let $\bR$ be the $J\times N$ matrix
with entries $R_{ji}$. Let $\bC\in\RealsP^J$
be a positive \emph{capacity} vector with components~$C_j$.
For each route $i$, \textit{packets} %(instead of flows, which is commonly
%used in the literature)
arrive as an independent Poisson process of rate $\lambda_i$. Packets arriving
on route $i$ require a unit amount of service, deterministically.

We denote the number of packets on route $i$ at time $t$ by $M_i(t)$,
and define the queue-size vector at time $t$ by $\bM
(t)=[M_i(t)]_{i=1}^N \in\IntegersP^N$.
Each packet gets service from the network at a rate determined according
to a bandwidth-sharing policy.
We also denote the total residual workload on route $i$ at time $t$ by
$W_i(t)$,
and let the vector of residual workload at time $t$ be $\bW(t) =
[W_i(t)]_{i=1}^N$.
Once a packet receives its total (unit) amount of
service, it departs the network.

We consider online, myopic bandwidth allocations. That is, the
bandwidth allocation
at time $t$ only depends on the queue-size vector $\bM(t)$. When there
are $m_i$
packets on route $i$, that is, if the vector of packets is $\mathbf{m}=
[m_i]_{i=1}^N$, let the
total bandwidth allocated to route $i$ be $\phi_i(\mathbf{m}) \in\RealsP
$. We consider
a processor-sharing policy, so that each packet on route $i$ is served
at rate
$\phi_i(\mathbf{m})/m_i$, if $m_i > 0$. If $m_i = 0$, let $\phi
_i(\mathbf{m}) = 0$.
If the \emph{bandwidth vector} $\bphi(\mathbf{m}) = [\phi_i(\mathbf{m}
)]_{i=1}^N$ satisfies the
capacity constraints
%
%e14 #&#
\begin{equation}
\label{eqfeasible} \bR\bphi(\mathbf{m}) \leq\bC\qquad\mbox{componentwise}
\end{equation}
for all $\mathbf{m}\in\IntegersP^{N}$, then, in light of Definition
\ref
{dfadmissible}, we say that
$\bphi(\cdot)$ is an \emph{admissible bandwidth allocation}.
A Markovian description of the system is given by a process
{$\bY(t)$} which contains the queue-size vector $\bM(t)$ along
with the residual workloads of the set of packets on each route.

Now, on average, $\lambda_i$ units of work arrive to route $i$ per unit
time. Therefore, in order for the Markov process $\bY(\cdot)$ to
be positive {(Harris)} recurrent, it is necessary that
%
%e15 #&#
\begin{eqnarray}
\label{eqcapfl} \bR\bolds{\lambda}& < & \bC\qquad\mbox{componentwise}.
\end{eqnarray}
All such $\bolds{\lambda}= [\lambda_i]_{i=1}^N \in\RealsP^N$
will be called
\textit{strictly admissible}, in the same spirit as strictly admissible
vectors for a switched network.
Similarly to the corresponding switched network, given $\bolds
{\lambda}\in\RealsP^N$,
we can define $\rho(\bolds{\lambda})$, the load induced by
$\bolds{\lambda}$, using~\eqref{eqdef-load-2},
as well as $\tilde{\rho}_j = (\sum_i R_{ji}\lambda_i )/C_j$.
Then by Lemma~\ref{lemsys-ind-load},
$\rho(\bolds{\lambda}) = \max_j \tilde{\rho}_j$, where $\tilde{\rho}_j$
can be interpreted
as the load induced by $\bolds{\lambda}$ on resource $j$.
%{[!I'VE CHANGED THE TERMS USED HERE TO USE ``ADMISSIBLE'', JUST TO
%LIMIT NOTATION A BIT.]}

\subsection*{Store-and-forward allocation (SFA) policy} We describe the
store-and-forward allocation policy that was first considered by
Massouli\'e and later analyzed in the thesis of Prouti\`ere \cite{PTh}.
Bonald and Prouti\`ere \cite{bonaldproutiere2} established that this
policy induces product-form
stationary distributions and is insensitive with respect to
phase-type distributions. This policy is shown to be \textit{insensitive}
for general service time distributions, including the deterministic service
considered here, by Zachary~\cite{zachary}. The relation between this policy,
the proportionally fair allocation, and multi-class queuing networks is discussed
in depth by Walton \cite{walton} and Kelly, Massouli{\'e} and Walton
\cite{KMW}. The insensitivity property implies
that the invariant measure of the process $\bM(t)$ only depends on the
parameters
$\bolds{\lambda}= [\lambda_i]_{i=1}^N \in\RealsP^N$, and no
other aspects of the stochastic description
of the system.

We first give an informal motivation for SFA. SFA is closely related
to quasi-reversible queuing networks. Consider a continuous-time
multi-class queuing \mbox{network} (without scheduling constraints) consisting
of processor sharing queues indexed by $j\in\{1,\ldots,J\}$ and job
types indexed by the routes $i\in\{1,\ldots,N\}$. Each route $i$ job
has a service requirement $R_{ji}$ at each queue $j$, and a fixed
service capacity $C_j$ is shared between jobs at the queue. Here each
job will sequentially visit all the queues (so-called
store-and-forward) and will visit each queue a fixed number of times.
If we assume that jobs on each route arrive as a Poisson process, then
the resulting queuing network will be stable for all strictly
admissible arrival rates. Moreover, each stationary queue will be
independent with a queue size that scales, with its load $\rho$, as
$\rho/(1-\rho)$. For further details, see Kelly \cite{Ke79}. So,
assuming each queue has equal load, the total number of jobs within the
network is of the order $J\rho/(1-\rho)$. In other words, these
networks have the stability and queue-size scaling that we require, but
do not obey the necessary scheduling constraints~\eqref{eqfeasible}.
However, these networks do produce an admissible schedule on average.
For this reason, we consider an SFA policy which, given the number of
jobs on each route, allocates the average rate with which jobs are
transferred through this multi-class network. Next, we describe this
policy (using notation similar to those used in \mbox{\cite{KMW,walton}}).

Given $\mathbf{m}\in\IntegersP^{N}$, define
\[
U(\mathbf{m}) = \biggl\{\tilde{\mathbf{m}} = \bigl(\tilde{m}_{ji}\dvtx  (j,
i) \in\cK\bigr) \in\IntegersP^{|\cK|}\dvtx  \sum_{j\dvtx  j \in i}
\tilde{m}_{ji} = m_i, \mbox{ for all } 1\leq i\leq N
\biggr\}.
\]
For $\bL\in\IntegersP^J$, we also define
\[
V(\bL) = \biggl\{\tilde{\mathbf{m}} = \bigl(\tilde{m}_{ji}\dvtx  (j, i) \in
\cK\bigr) \in\IntegersP^{|\cK|}\dvtx  \sum_{i\dvtx  i \ni j}
\tilde{m}_{ji} = L_j,\mbox{ for all } 1\leq j\leq J \biggr
\}.
\]
Here, by notation $j\in i$ (and $i \ni j$) we mean $R_{ji} > 0$.
%The notation $i \ni j$ is used
%when we consider a collection of $i$ satisfying this condition for a
%given $j$. }
For each
$\tilde{\mathbf{m}} \in U(\mathbf{m})$, we exploit notation somewhat
and define
$\tilde{m}_j = \sum_{i\dvtx  j \in i} \tilde{m}_{ji}$, for
all $j\leq J$. Also define
\[
\pmatrix{\tilde{m}_j \vspace*{2pt}\cr\tilde{m}_{ji}\dvtx  i \ni j} =
\frac{\tilde{m}_j !}{\prod_{i\dvtx  j \in i} (\tilde{m}_{ji}!)}.
\]
%
%In the above, by $i \ni j$ we mean that $R_{ji} > 0$; the notation $i
%when we consider a collection of $i$ satisfying this condition for a
%given $j$.
For
$\mathbf{m}\in\IntegersP^{N}$, we define $\Phi(\mathbf{m})$ as
%
%e16 #&#
\begin{equation}
\label{dfsfa1} \Phi(\mathbf{m}) = \sum_{\tilde{\mathbf{m}} \in U(\mathbf
{m})} \prod
_{j \in J} \biggl(\pmatrix{\tilde{m}_j \vspace*{2pt}\cr
\tilde{m}_{ji}\dvtx  i \ni j}\prod_{i\dvtx  j\in i} \biggl(
\frac{R_{ji}}{C_j} \biggr)^{\tilde{m}_{ji}} \biggr).
\end{equation}
We shall define $\Phi(\mathbf{m}) = 0$ if any of the components of
$\mathbf{m}$
is negative.
The store-and-forward allocation (SFA) assigns rates according to the function
$\bphi\dvtx  \IntegersP^{N} \to\RealsP^{N}$, so that for any
$\mathbf{m}\in\IntegersP^{N}$, $\bphi(\mathbf{m}) = (\phi_i(\mathbf
{m}))_{i=1}^N$, with
%
%e17 #&#
\begin{equation}
\label{dfsfa2} \phi_i(\mathbf{m}) = \frac{\Phi(\mathbf{m}- \be_i)}{\Phi
(\mathbf{m})},
\end{equation}
where, recalling that $\mathbf{m}-\be_i$ is the same as $\mathbf{m}$ at
all but
the $i$th component,
its $i$th component equals $m_i - 1$. The bandwidth allocation {$\bphi
(\mathbf{m})$} is the stationary throughput of jobs on the routes of a
multi-class queuing network (described above), conditional on there
being $\mathbf{m}$ jobs on each route.

A priori it is not clear if the above
described bandwidth allocation is even admissible, that is, satisfies
\eqref{eqfeasible}. {This can be argued as follows. The
$\bphi(\mathbf{m})$ can be related to the stationary throughput of a
multi-class network with a finite number of jobs, $\mathbf{m}$, on each route.
Under this scenario (due to finite number of jobs), each queue must be
stable. Therefore, the load on each queue,
$\bR\bphi(\mathbf{m})$, must be less than the overall system capacity
$\bC
$. That is, the allocation is admissible.
The precise argument along these lines is provided in, for example,
\cite{KMW}, Corollary 2 and \cite{walton}, Lemma 4.1.
}

The SFA induces a product-form invariant distribution for the number of
packets waiting
in the bandwidth-sharing network and is insensitive. %\cite{zachary}.
We summarize this in the following result. %known result.
%
%th4.1 #&#
\begin{theorem}\label{thmSFA1}%[\cite{bonaldproutiere1}, {
Consider a bandwidth-sharing network with $\bR\bolds{\lambda}<
\bC$. Under the
SFA policy described above, the Markov process $\bY(t)$ is positive
{(Harris)} recurrent,
and $\bM(t)$ has a unique stationary probability distribution ${\bolds
{\pi}}$
given by
%
%e18 #&#
\begin{equation}
\label{eqsfameasure} {\bolds{\pi}}(\mathbf{m}) = \frac{\Phi(\mathbf
{m})}{\Phi}\prod
_{i=1}^N \lambda_i^{m_i}\qquad\mbox{for all } \mathbf{m}\in\IntegersP^{N},
\end{equation}
where
%
%e19 #&#
\begin{equation}
\label{eqnormalizeSFA} \Phi= \prod_{j=1}^J \biggl(
\frac{C_j}{C_j-\sum_{i\dvtx  i\ni j} R_{ji}
\lambda_i} \biggr)
\end{equation}
is a normalizing factor.
Furthermore, the steady-state residual workload of packets waiting in
the network
can be characterized as follows.
First, the steady-state distribution of the residual workload of a
packet is
independent from ${\bolds{\pi}}$.
Second, in steady state, conditioned on the number of packets on each
route of the network,
the residual workload of each packet is uniformly distributed on $[0, 1]$,
and is independent from the residual workloads of other packets.
%%%\iffalse
%%%{Furthermore, in steady state, the residual workload of each packet
%in
%%%the network
%%%is uniformly distributed on $[0, 1]$
%%%and} {is conditionally independent from the residual workloads of
%other
%%%packets, when we condition on the number of packets on each route of
%%%the network.}\fi
\end{theorem}
Note that statements similar to Theorem~\ref{thmSFA1} have appeared
in other works, for example,
\cite{bonaldproutiere1}, \cite{walton}, Proposition 4.2, and \cite{KMW}.
Theorem~\ref{thmSFA1} is a summary of these statements, and for completeness,
it is proved in Appendix~\ref{apdxpsn}.

The following property of the stationary distribution ${\bolds{\pi}}$ described
in Theorem~\ref{thmSFA1} will be useful.
%
%pr4.2 #&#
\begin{proposition}\label{propSFA2}
Consider the setup of Theorem~\ref{thmSFA1}, and let ${\bolds{\pi}}$ be
described by \eqref{eqsfameasure}.
Define a measure $\tilde{{\bolds{\pi}}}$ on $\IntegersP^{|\cK|}$ as follows:
for $\tilde{\mathbf{m}} \in\IntegersP^{|\cK|}$,
%
%e20 #&#
\begin{eqnarray}
\label{eqmcmeasure} \tilde{{\bolds{\pi}}}(\tilde{\mathbf{m}}) & =& \frac
{1}{\Phi}
\prod_{j=1}^J \biggl(\pmatrix{\tilde{m}_j
\vspace*{2pt}\cr\tilde{m}_{ji}\dvtx  i \ni j} \prod_{i\dvtx  j\in i}
\biggl(\frac{R_{ji} \lambda_i}{C_j} \biggr)^{\tilde{m}_{ji}} \biggr).
\end{eqnarray}
%
%{DELETE?: In the above, we have used the notation $\tilde{m}_j =
%i \ni j} \tilde{m}_{ji}$, for any $j$.}
Then, for any $L \in\IntegersP$,
%
%e21 #&#
\begin{eqnarray}
\label{eqsfaeqmc} {\bolds{\pi}} \Biggl( \Biggl\{\mathbf{m}\dvtx  \sum
_{i=1}^{N} m_i = L \Biggr\} \Biggr) & =&
\tilde{{\bolds{\pi}}} \Biggl( \Biggl\{\tilde{\mathbf{m}}\dvtx  \sum
_{j=1}^{J} \tilde{m}_j = L \Biggr\}
\Biggr).
\end{eqnarray}
\end{proposition}
We relate the distribution $\tilde{{\bolds{\pi}}}$ to the stationary
distribution of
an insensitive multi-class queuing network with a product-form
stationary distribution and geometrically distributed queue sizes.
%
%pr4.3 #&#
\begin{proposition}\label{propMC}
Consider the distribution $\tilde{{\bolds{\pi}}}$ defined in \eqref
{eqmcmeasure}. Then,
for any $\bL= (L_1,\ldots, L_{J}) \in\IntegersP^{J}$,
%
%e22 #&#
\begin{eqnarray}\label{eqmcindep}
\qquad \tilde{{\bolds{\pi}}} (\tilde{m}_1 = L_1,\ldots, \tilde{m}_{J} = L_{J} ) & %\stackrel{(a)}{=}
=& \sum
_{(\tilde{m}_{ji}) \in V(\bL)} \tilde{{\bolds{\pi}}}\bigl((\tilde{m}_{ji})
\bigr)
 = \prod_{j=1}^{J} \tilde{
\rho}_j^{L_j} (1-\tilde{\rho}_j),
\end{eqnarray}
where $\tilde{\rho}_j = (\sum_{i\dvtx  i\ni j} R_{ji} \lambda_j )/C_j$.
%{PERHAPS WE DO NOT NEED TO SAY THIS?: In the above, (a) is effectively
%the definition of $\tilde{{\bolds{\pi}}}$ when projected
%to the $J$-dimensional space of $(\tilde{m}_1,\ldots, \tilde{m}_{J})$.}
\end{proposition}

Using Theorem~\ref{thmSFA1} and Propositions~\ref{propSFA2} and~\ref{propMC},
we can compute the expected value and the probability tail exponent of
the steady-state total residual workload in the system. Recall that the
total residual workload in the system
at time $t$ is given by $\sum_{i=1}^N W_i(t)$.
%
%pr4.4 #&#
\begin{proposition}\label{propworkload}
Consider a bandwidth-sharing network with $\bR\bolds{\lambda}<
\bC$,
operating under the SFA policy.
Denote the load induced by $\bolds{\lambda}$ to be $\rho= \rho
(\bolds{\lambda}) (< 1)$,
and for each $j$, let $\tilde{\rho}_j = (\sum_{i} R_{ji}\lambda
_i )/C_j$.
Then $\bW(\cdot)$ has a unique stationary probability distribution.
With respect to this stationary distribution, the following properties hold:
\begin{longlist}[(ii)]
\item[(i)] The expected total residual workload is given by
%
%e23 #&#
\begin{equation}
\label{eqmean-workload} \bbE\Biggl[\sum_{i=1}^N
W_i \Biggr] = \frac{1}{2}\sum_{j=1}^J
\frac{\tilde{\rho}_j}{1-\tilde{\rho}_j}.
\end{equation}
\item[(ii)] The distribution of the total residual workload has an
exponential tail
with exponent given by
%
%e24 #&#
\begin{equation}
\label{eqtail-workload} \lim_{L \to\infty} \frac{1}{L} \log\Prob\Biggl(
\sum_{i=1}^N W_i \geq L
\Biggr) = -\theta^*,
\end{equation}
where $\theta^*$ is the unique \emph{positive} solution of the equation
$\rho(e^{\theta} - 1) = \theta$.
\end{longlist}
\end{proposition}
\section{Main result: A policy and its performance}\label{secmain}

In this section, we describe an online scheduling policy and
quantify its performance in terms of explicit, closed-form bounds on the
stationary distribution of the induced queue sizes. Section~\ref{ssecpolicy} describes\vadjust{\goodbreak} the policy for a generic switched
network and provides the statement of the main result. Section~\ref{sseciq}
discusses its implications. Specifically, it discusses (a) the optimality
of the policy for a large class of switched networks with respect to exponential
tail bounds, and (b) the optimality of the policy for a class of switched
networks, including input-queued switches, with respect to the average
total queue size. Section~\ref{ssecproof} proves the main result
stated in
Section~\ref{ssecpolicy}.

%s5.1 #&#
\subsection{A policy for switched networks}\label{ssecpolicy}

The basic idea behind the policy, to be described in detail shortly,
is as follows. Given a switched network, denoted by~$\bSN$, with
constraint set $\sS$
and $N$ queues, let $\hull{\sS}$ have rank $J$ and representation
[cf. \eqref{eqrankeq}]
\[
\hull{\sS}= \bigl\{ \mathbf{x}\in[0,1]^N \dvtx  \bR\mathbf{x}\leq\bC
\bigr\}, \qquad\bR\in\RealsP^{N \times J},  \bC\in\RealsP^{J}.
\]
Now consider a virtual bandwidth-sharing network, denoted by $\bBN$,
with $N$ routes corresponding to each of these $N$ queues. The
resource--route relation is determined precisely by the matrix $\bR$,
and the $J$
resources have capacities given by $\bC$. Both networks, $\bSN$ and
$\bBN$,
are fed identical arrivals. That is, whenever a packet arrives to queue
$i$ in $\bSN$, a packet is added to route $i$ in $\bBN$ at the same time.
The main question is that of determining a scheduling policy for $\bSN$;
this will be derived from $\bBN$. Specifically, $\bBN$ will operate
under the insensitive SFA policy described in Section~\ref{secinsensitive}.
By Theorem~\ref{thmSFA1} and Propositions~\ref{propSFA2} and
\ref{propMC}, this will induce a desirable stationary distribution of
queue sizes in $\bBN$. Therefore, if we could use the rate allocation
of $\bBN$, that is, the SFA policy, directly in $\bSN$, it would give us
a desired performance in terms of the stationary distribution of the
induced queue sizes.
Now the rate allocation in $\bBN$ is such that the instantaneous rate
is always inside $\hull{\sS}$. However, it could change all the time
and need not utilize points of $\sS$ as rates. In contrast,
in $\bSN$ we require that the rate allocation can change only once
per discrete time slot and it must always employ one of the generators
of $\hull{\sS}$, that is, a schedule from $\sS$. The key to our policy
is an effective way to emulate the rate allocation of $\bBN$ under SFA
(or for that matter, any admissible bandwidth allocation) by utilizing schedules
from $\sS$ in an online manner and with the discrete-time constraint.
We will
see {shortly} that this emulation policy relies on $\sS$ being
monotone; cf. Assumption~\ref{assmonotone}.

To that end, we describe this emulation policy. Let us start by introducing
some useful notation. Let $\bA(\cdot) = (A_i(\cdot) )$ be
the vector of exogenous,
independent Poisson processes according to which unit-sized packets
arrive to both
$\bBN$ and $\bSN$, simultaneously. Recall that $A_i(\cdot)$ is a
Poisson process
with rate $\lambda_i$. Let $\bM(t) = (M_i(t) )$ denote the
vector of
numbers of packets waiting on the $N$ routes in $\bBN$ at time $t \geq
0$. In
$\bBN$, the services are allocated according to the SFA policy
described in Section~\ref{secinsensitive}.
Let $\Lambda^{\mathrm{SFA}}(\cdot) = (\Lambda^{\mathrm{SFA}}_i(\cdot)
)\in\RealsP^N$
denote the cumulative amount of service allocated to the $N$ routes
in $\bBN$ under the SFA policy: $\Lambda^{\mathrm{SFA}}_i(t)$ denotes
the total
amount of service allocated to all packets on route $i$ during the
interval $[0,t]$,
for $t \geq0$, with $\Lambda^{\mathrm{SFA}}_i(0) = 0$ for $1\leq i\leq
N$. By definition,
all components of $\Lambda^{\mathrm{SFA}}(\cdot)$ are nondecreasing
and Lipschitz
continuous. Furthermore,
$(\Lambda^{\mathrm{SFA}}(t+s) - \Lambda^{\mathrm{SFA}}(t))/s \in\hull
{\sS}$ for any $t \geq0$ and
$s > 0$.
Recall that the \mbox{(right-)derivative} of $\Lambda^{\mathrm{SFA}}(\cdot)$
is determined by
$\bM(\cdot)$ through the function $\bphi(\cdot)$ as defined in
\eqref{dfsfa2}.
%$\pLSFA(\cdot)$ is determined by $\bM(\cdot)$.

Now we describe the scheduling policy for $\bSN$ that will rely on
$\Lambda^{\mathrm{SFA}}(\cdot)$.
Let $\bB(\tau) = (B_i(\tau) )$ denote the cumulative amount
of service
allocated in $\bSN$ by the scheduling policy up to time slot $\tau
\geq0$, with
$\bB(0) = \bZero$.
%
%Policy will also maintain {\em virtual} cumulative amount of service $
%and $\tbB(\tau) \leq\bB(\tau)$ componentwise for all $\tau\geq0$.
%
The scheduling policy determines how $\bB(\cdot)$ is updated. Let
$\bQ(\tau) = (Q_i(\tau) )$ be the queue sizes measured
at the end of time slot $\tau$. Let service be provided according to
the scheduling
policy instantly at the beginning of a time slot. Thus, the scheduling
policy decides
the schedule $d\bB(\tau) = \bB(\tau+1) - \bB(\tau) \in\sS$
% (as well as $d\tbB(\tau) = \tbB(\tau+1) - \tbB(\tau)$)
at the very beginning of time slot $\tau+1$. This decision is made as
follows. Let $\bD(\tau) = \Lambda^{\mathrm{SFA}}(\tau) - \bB(\tau)$.
We will see shortly that under our policy, $\bD(\tau)$ is always
nonnegative.
This fact will be useful at various places, and in particular, for
bounding the discrepancy
between the continuous-time policy SFA and its discrete-time emulation.
Let $\rho(\bD(\tau))$ be
the optimal objective value in the optimization problem $\PRIMAL(\bD
(\tau))$ defined in
\eqref{defprimal}. In particular, there exists a nonnegative
combination of
schedules in $\sS$ such that
%
%e25 #&#
\begin{equation}
\label{eqcvx0} \sum_{\bolds{\sigma}\in\sS} \tilde{
\alpha}_{\bolds{\sigma}} \bolds{\sigma}\geq\bD(\tau) \quad\mbox{and}\quad\sum
_{\bolds{\sigma}\in\sS} \tilde{\alpha}_{\bolds{\sigma}
} = \rho\bigl(\bD(
\tau)\bigr).
\end{equation}
We claim that in fact, we can find nonnegative numbers $\alpha
_{\bolds{\sigma}}$, $\bolds{\sigma}\in\sS$,
such that
%
%e26 #&#
\begin{equation}
\label{eqcvx} \sum_{\bolds{\sigma}\in\sS} \alpha_{\bolds{\sigma}}
\bolds{\sigma}= \bD(\tau) \quad\mbox{and}\quad\sum_{\bolds{\sigma}\in\sS}
\alpha_{\bolds{\sigma}} = \rho\bigl(\bD(\tau)\bigr).
\end{equation}
This is formalized in the following lemma.
%
%le5.1 #&#
\begin{lemma}\label{LEMMONOTONE}
Let $\bD\in\RealsP^N$ be a nonnegative vector.
Consider the static planning problem $\PRIMAL(\bD)$ defined in \eqref
{defprimal}.
Let the optimal objective value to $\PRIMAL(\bD)$ be $\rho(\bD)$.
Then there exists $\alpha_{\bolds{\sigma}} \geq0$, $\bolds{\sigma}\in\sS$,
such that \eqref{eqcvx} holds.
\end{lemma}
The proof of the lemma relies on Assumption~\ref{assmonotone},
and is provided in the \hyperref[app]{Appendix}.

There could be many possible nonnegative combinations of $\bD(\tau)$
satisfying~\eqref{eqcvx}. If there exist nonnegative numbers $\alpha_{\bolds
{\sigma}
}$, $\bolds{\sigma}\in\sS$,
satisfying \eqref{eqcvx} with $\alpha_{\bolds{\sigma}'} \geq1$ for some
$\bolds{\sigma}' \in\sS$, then choose $\bolds{\sigma}'$ as the
schedule: set
$d\bB(\tau) = \bolds{\sigma}'$. If no such decomposition exists for $\bD
(\tau)$,
then set $d\bB(\tau) = \tilde{\bolds{\sigma}}$, where $\tilde{\bolds
{\sigma}}$
is a solution (ties broken arbitrarily) of
%
%e27 #&#
\begin{equation}
\mbox{maximize}\qquad  \sum_i
\sigma_i \quad\mbox{over}\quad\bolds{\sigma}\in\sS, \bolds{\sigma}
\leq\bD(\tau).
\end{equation}
Here first observe that for all time $\tau$,
$d\bB(\tau) \leq\bD(\tau)$, so $\bD(\tau) \geq\bZero$.
Hence, $\bZero$ is a feasible solution for the above problem, as
$\bZero\in\sS$.

The above is a complete description of the scheduling policy. Observe
that it is an online policy, as the virtual network $\bBN$ can be simulated
in an online manner, and, given this, the scheduling decision in $\bSN
$ relies
only on the history of $\bBN$~and~$\bSN$.
The following result quantifies the performance of the policy.

%%%\iffalse
%%%Therefore, $\bD(\tau)/\rho(\bD(\tau)) \in\hull{\sS}$. Since
%%%$\hull{\sS}$
%%%is an $N$ dimensional bounded convex set, it follows by Carath
%%%theorem that we can have decomposition satisfying \eqref{eqcvx}
%%%such that $|\{\bolds{\sigma}\in\sS\dvtx  \alpha_{\bolds{\sigma}} > 0\}|
%%%Consider one such convex decomposition of $\bD(\tau)$. If
%%%\fi

%
%th5.2 #&#
\begin{theorem}\label{thmmain}
Given a strictly admissible arrival rate vector $\bolds{\lambda
}$, with $\rho= \rho(\bolds{\lambda}) < 1$, under
the policy described above, the switched network $\bSN$ is positive recurrent
and has a unique stationary distribution.
Let $\tilde{\rho}_j=\break  (\sum_{i} R_{ji} \lambda_i )/C_j$, $j =
1, 2, \ldots, J$
be the same as in Proposition~\ref{propworkload}. With respect to this
stationary
distribution, the following properties hold:
\begin{longlist}[(2)]
\item[(1)] The expected total queue size is bounded as
%
%e28 #&#
\begin{eqnarray}
\label{eqthmmain1} \bbE\Biggl[ \sum_{i=1}^N
Q_i \Biggr] & \leq&\frac{1}{2} \Biggl(\sum
_{j=1}^{J} \frac{\tilde{\rho}_j}{1-\tilde{\rho}_j} \Biggr) + K (N+2),
\end{eqnarray}
where %$\tilde{\rho}_j = \big(\sum_{i} R_{ji} \lambda_i\big)/C_j$ are as
%in Proposition~\ref{propMC},
$K = \max_{\bolds{\sigma}\in\sS} (\sum_i \sigma_i )$.

\item[(2)] The distribution of the total queue size has an exponential
tail with exponent given by
%
%e29 #&#
\begin{eqnarray}
\label{eqthmmain2} \lim_{L \to\infty} \frac{1}{L} \log\Prob\Biggl(
\sum_{i=1}^N Q_i \geq L
\Biggr) & =& - \theta^*,
\end{eqnarray}
where $\theta^*$ is the unique positive solution of the equation $\rho
(e^{\theta} - 1) = \theta$.
\end{longlist}
\end{theorem}

%s5.2 #&#
\subsection{Optimality of the policy}\label{sseciq}

This section establishes the optimality of our policy for
input-queued switches, both with respect to expected total queue-size scaling and tail exponent.
{General conditions under which our policy is optimal
with respect to tail exponent are also provided.}
%The policy produces an optimal
%tail exponent for any single-hop switched network. %However, the
%average
%queue-size scaling, which depends on the {\em rank} of $\sS$,
%is desirable for systems when {\em rank} of $\sS$ is
%smaller than the number of queues $N$; indeed,
%input-queued switch is such a case.

\subsubsection*{Scaling of queue sizes}\label{sssecformal}
We start by formalizing what we mean by the optimality of expected
queue sizes
and of their tail exponents. We consider
policies under which there is a well-defined {limiting} stationary distribution
of the queue sizes for all $\bolds{\lambda}$ such that $\rho
(\bolds{\lambda}) < 1$.
Note that the class of policies is not empty;
indeed, the maximum weight policy and our policy are members of
this class. With some abuse of notation, let
${\bolds{\pi}}$ denote the stationary distribution of the queue-size vector
under the policy of interest. We are interested in two quantities:
\begin{longlist}[(2)]
\item[(1)] \textit{Expected total queue size}. Let $\widebar{Q}$ be the expected
total queue size under the stationary distribution ${\bolds{\pi}}$,
defined by
\[
\widebar{Q} = \bbE_{{\bolds{\pi}}} \biggl[\sum_i
Q_i \biggr].\vadjust{\goodbreak}
\]
Note that by ergodicity, the time average of the total queue size
and the expected total queue size under ${\bolds{\pi}}$ are the same quantity.

\item[(2)] \textit{Tail exponent}. Let $\beta_L(Q), \beta_U(Q) \in
[-\infty, 0]$
be the lower and upper limits of the tail exponent of the total queue size
under ${\bolds{\pi}}$ (possibly $-\infty$ or $0$), respectively,
defined by
%
%e30 #&#
%e31 #&#
\begin{eqnarray}
\beta_L(Q) & =& \lim\inf_{\ell\to\infty} \frac{1}{\ell}
\log\Prob_{{\bolds{\pi}}} \biggl(\sum_i
Q_i \geq\ell\biggr)\quad\mbox{and} \label{eqlexp}
\\
\beta_U(Q) & =& \lim\sup_{\ell\to\infty}
\frac
{1}{\ell} \log\Prob_{{\bolds{\pi}}} \biggl(\sum
_i Q_i \geq\ell\biggr). \label{equexp}
\end{eqnarray}
If $\beta_L(Q) = \beta_U(Q)$, then we denote this common value by
$\beta(Q)$.
\end{longlist}
We are interested in policies that can achieve minimal $\widebar{Q}$ and
$\beta(Q)$. For tractability, we focus on
\textit{scalings} of these quantities with respect to $\sS$
(equivalently, $N$)
and $\rho(\bolds{\lambda})$, as $1/(1-\rho(\bolds{\lambda
}))$ and $N$ increase.
For different $\bolds{\lambda}'$ and $\bolds{\lambda}$, it
is possible that
$\rho(\bolds{\lambda}) = \rho(\bolds{\lambda}')$, but the
scaling of $\widebar{Q}$, for example,
could be wildly different. For this reason, we consider the worst possible
dependence on $1/(1-\rho)$ and $N$ among all $\bolds{\lambda}$
with $\rho(\bolds{\lambda}) = \rho$.
%With this understanding, we shall call a policy optimal with respect
%to scaling
%if its dependence on $1/(1-\rho)$ and $N$ is the best possible up to a
%universal
%constant, as $1/(1-\rho)$ and $N$ are increasing.

Note that we are considering scalings with respect to two quantities,
$\rho$ and $N$, and we are interested in two limiting regimes, $\rho
\to1$
and $N \to\infty$. The optimality of queue-size scaling stated here
is with
respect to the order of limits $\rho\to1$ and then $N\to\infty$. As
noted in \cite{STZopen},
taking the limits in different orders could potentially result in different
limiting behaviors of the object of interest, for example, $\widebar{Q}$.
For further discussion, see Section~\ref{secconc}. It should be noted,
however, that
whenever the tail exponent is optimal, this optimality holds for \textit{any} $\rho$ and $N$.

\subsubsection*{Optimality of the tail exponent}\label{sssecgeneral}
%{[THIS WHOLE PARAGRAPH IS REVISED.]}
Here we establish sufficient conditions
under which our policy is optimal with respect to tail exponent.
%the optimality of the tail exponent
%for input-queued switches under our policy.
First, we present a universal lower bound on the tail exponent,
for a general single-hop switched network under any policy.
We then provide a condition under which this lower bound matches
the tail exponent under our policy.
This condition is satisfied by both input-queued switches
and the independent-set model of wireless networks.
%This lower bound is then specialized to the context of input-queued
%switches,
%and compared against the tail exponent under our policy.

Consider any policy under which there exists
a well-defined limiting stationary distribution of the queue sizes
for all $\bolds{\lambda}$ such that $\rho(\bolds{\lambda})
< 1$.
Let ${\bolds{\pi}}_0$ denote the stationary distribution of queue sizes
under this policy. The following lemma establishes
a universal lower bound on the tail exponent.

%
%le5.3 #&#
\begin{lemma}\label{lemtail}
Consider a switched network as described in Theorem~\ref{thmmain},
with scheduling set $\sS$ and admissible region $\{\mathbf{x}\in[0,
1]^N \dvtx  \bR\mathbf{x}\leq\bC\}$.
Let ${\bolds{\pi}}_0$ and $\bolds{\lambda}$ be as described.
For each $j$, let $\tilde{\rho}_j = \sum_{i=1}^N R_{ji}\lambda_i / C_j$
be defined as in Theorem~\ref{thmmain}.
Then under ${\bolds{\pi}}_0$,
%
%e32 #&#
\begin{equation}
\label{eqtail-lower1} \liminf_{L \to\infty} \frac{1}{L} \log
\Prob_{{\bolds{\pi}}_0} \biggl(\sum_{i}
Q_i \geq L \biggr) \geq- \min_{j = 1, 2, \ldots, J}
\theta^*_j,
\end{equation}
where, for each $j \in\{1, 2, \ldots, J\}$, $\theta^*_j$
is the unique positive solution of the equation
\[
\sum_{i=1}^N \lambda_i
\bigl(e^{R_{ji}\theta} - 1 \bigr) = \theta.
\]
\end{lemma}
\begin{pf}
Consider a fixed $j \in\{1, 2, \ldots, J\}$.
Without loss of generality, we assume that $C_j = 1$, by properly normalizing
the inequality $(\bR\mathbf{x})_j \leq C_j$. In this case, $R_{ji}
\leq1$ for all $i$,
since for each $i \in\{1, 2, \ldots, N\}$,
$\be_i \in\sS\subset\hull{\sS}$, and satisfies the constraint
$(\bR\be_i)_j = R_{ji} \leq C_j = 1$.

Now consider the following single-server queuing system.
The arrival process is given by the sum $\sum_{i = 1}^N R_{ji}
A_i(\cdot)$,
so that arrivals across time slots are independent,
and that in each time slot, the amount of work that arrives is $\sum
_{i=1}^N R_{ji} a_i$,
where $a_i$ is an independent Poisson random variable with mean
$\lambda_i$, for each $i$.
Note that the arriving amount in a single time slot does not have to be
integral.
Note also that $\sum_{i=1}^N R_{ji}\lambda_i = \tilde{\rho}_j < 1$,
since $\rho(\bolds{\lambda}) = \max_j \tilde{\rho}_j < 1$.
In each time slot, a unit amount of service is allocated to the total
workload in the system.
Then, for this system, the workload process $W(\cdot)$ satisfies
\[
W(\tau+1) = \bigl[W(\tau) - 1 \bigr]^+ + \sum_{i=1}^N
R_{ji}a_i(\tau),
\]
where $a_i(\tau)$ is the number of arrivals to queue $i$ in the
original system in time slot $\tau$.
We make two observations for this system.
First, $W(\cdot)$ is stochastically dominated by $\sum_{i=1}^N R_{ji}
Q_i(\cdot)$,
where $Q_i(\cdot)$ is the size of queue $i$ in the original system,
under any online scheduling\vspace*{1pt} policy.
This is because for all schedules $\bolds{\sigma}\in\sS$, $\bolds{\sigma
}$ satisfies
$\bR\bolds{\sigma}\leq\bC$,
and hence $\sum_{i=1}^N R_{ji}\sigma_i \leq C_j = 1$ for every $\bolds
{\sigma}
\in\sS$.
Second, since $R_{ji} \leq1$ for all $i$,
$\sum_{i=1}^N R_{ji} Q_i(\cdot)$ is stochastically dominated by $\sum
_{i=1}^N Q_i(\cdot)$.
Thus we have
\[
\liminf_{L \to\infty} \frac{1}{L} \log\Prob_{{\bolds{\pi}}_0}
\biggl(\sum_{i} Q_i \geq L \biggr) \geq
\liminf_{L \to\infty} \frac{1}{L} \log\Prob\bigl(W(\infty) \geq
L \bigr).
\]
We now show that
\[
\liminf_{L \to\infty} \frac{1}{L} \log\Prob\bigl(W(\infty)
\geq L \bigr) \geq-\theta^*_j, %- \frac{2(1-\tilde{\rho}_j)}{
\]
where $\theta^*_j$ is the unique positive solution of the equation
\[
\sum_{i=1}^N \lambda_i
\bigl(e^{R_{ji}\theta} - 1 \bigr) = \theta.
\]
Consider the log-moment generating function (log-MGF)
of $\sum_{i=1}^N R_{ji}a_i$, the arriving amount in one time slot.
Since $a_i$ is a Poisson random variable with mean $\lambda_i$ for
each $i$,
its moment generating function is given by
\[
f(\theta) = \exp\Biggl(\sum_{i=1}^N
\lambda_i\bigl(e^{R_{ji}\theta
}-1\bigr) \Biggr).
\]
Hence the log-MGF is
\[
\log f(\theta) = \sum_{i=1}^N
\lambda_i\bigl(e^{R_{ji}\theta}-1\bigr).
\]
By Theorem 1.4 of \cite{bigQ},
\[
\lim_{L \to\infty} \frac{1}{L} \log\Prob\bigl(W(\infty) \geq L
\bigr) = -\theta^*_j,
\]
where $\theta^*_j = \sup\{\theta> 0\dvtx  \log f(\theta) < \theta\}$.
Since $\log f(\theta) - \theta$ is strictly convex, $\theta^*_j$
satisfies
\[
\sum_{i=1}^N \lambda_i
\bigl(e^{R_{ji}\theta^*_j} - 1 \bigr) = \theta^*_j.
\]
$j \in\{1, 2, \ldots, J\}$ is arbitrary, so
\[
\liminf_{L \to\infty} \frac{1}{L} \log\Prob_{{\bolds{\pi}}_0}
\biggl(\sum_{i} Q_i \geq L \biggr)
\geq- \min_{j = 1, 2, \ldots, J} \theta^*_j.
\]
\old{
Thus, if we find some $\tilde{\theta}$ such that $\log f(\tilde{\theta})
\geq\tilde{\theta}$,
then
\[
\liminf_{L \to\infty} \frac{1}{L} \log\Prob\bigl(W(\infty)
\geq L \bigr) \geq-\tilde{\theta}.
\]
Let $\tilde{\theta} = \frac{2(1-\tilde{\rho}_j)}{\sum_{i=1}^N \lambda
_i R_{ji}^2}$.
Then,
\begin{eqnarray*}
\log f(\tilde{\theta}) &=& \sum_{i=1}^N
\lambda_i\bigl(e^{R_{ji}\tilde{\theta}}-1\bigr) \geq\sum
_{i=1}^N \lambda_i\biggl(R_{ji}
\tilde{\theta} + \frac{R_{ji}^2\tilde{\theta}^2}{2}\biggr)
\\
& = & \tilde{\rho}_j\tilde{\theta} + \frac{\sum_{i=1}^N \lambda_i
R_{ji}^2}{2} \tilde{
\theta}^2
\\
& = & \tilde{\rho}_j \frac{2(1-\tilde{\rho}_j)}{\sum_{i=1}^N \lambda_i
R_{ji}^2} + \frac{\sum_{i=1}^N \lambda_i R_{ji}^2}{2}
\biggl(\frac{2(1-\tilde{\rho}_j)}{\sum_{i=1}^N \lambda_i
R_{ji}^2} \biggr)^2
\\
& = & \frac{2\tilde{\rho}_j(1-\tilde{\rho}_j)}{\sum_{i=1}^N \lambda_i
R_{ji}^2} + \frac{2(1-\tilde{\rho}_j)^2}{\sum_{i=1}^N \lambda_i R_{ji}^2}
\\
& = & \frac{2(1-\tilde{\rho}_j)}{\sum_{i=1}^N \lambda_i R_{ji}^2} =
\tilde{\theta}.
\end{eqnarray*}
Thus, we have shown that
\[
\liminf_{L \to\infty} \frac{1}{L} \log\Prob\bigl(W(\infty)
\geq L \bigr) \geq-\frac{2(1-\tilde{\rho}_j)}{\sum_{i=1}^N \lambda_i R_{ji}^2}.
\]
Since this holds for an arbitrary $j \in\{1, 2, \ldots, J\}$,
\[
\liminf_{L \to\infty} \frac{1}{L} \log\Prob\bigl(W(\infty)
\geq L \bigr) \geq\max_{j = 1, 2, \ldots, J} -\frac{2(1-\tilde{\rho
}_j)}{\sum_{i=1}^N \lambda_i R_{ji}^2}.
\]
}\upqed
\end{pf}
For general switched networks, the lower bound above need not match
the tail exponent
achieved under our policy [cf. \eqref{eqthmmain2}]. However, for a
wide class of switched networks, these two quantities are equal.
The following corollary of Lemma~\ref{lemtail} is immediate.
%
%co5.4 #&#
\begin{corollary}\label{corsuff-opt-tail-expo}
Consider a switched network as described in Lemma~\ref{lemtail},
with scheduling set $\sS$ and admissible region $\{\mathbf{x}\in[0,
1]^N \dvtx  \bR\mathbf{x}\leq\bOne\}$.
If for all $j$ and $i$, $R_{ji} \in\{0, 1\}$,
then our policy achieves optimal tail exponent,
for any strictly admissible arrival-rate vector $\bolds{\lambda}$.
\end{corollary}
\begin{pf}
Let $\bolds{\lambda}\in\RealsP^N$ be strictly admissible, that
is, $\bR\bolds{\lambda}< \bOne$.
Let $\tilde{\rho}_j = \sum_{i} R_{ji} \lambda_i$ for each $j$,
and let $\rho= \rho(\bolds{\lambda})$ be the system load
induced by $\bolds{\lambda}$.
Consider the $\theta^*_j$ in Lemma~\ref{lemtail}.
When $R_{ji} \in\{0, 1\}$ for all $j$, and $i$, $\theta^*_j$
is the unique positive solution of the equation
\[
\tilde{\rho}_j\bigl(e^{\theta} - 1\bigr) = \theta
\]
for each $j$. Using the relation $\rho= \max_j \tilde{\rho}_j$,
we see that $\min_{j} \theta^*_j$ is the unique positive solution of
the equation
\[
\rho\bigl(e^{\theta} - 1\bigr) = \theta.
\]
Comparing this with equation \eqref{eqthmmain2} of Theorem~\ref{thmmain},
we see that our policy achieves the optimal tail exponent.
\end{pf}
Consider an $n\times n$ input-queued switch,
defined in Section~\ref{sseciqintro}, and with the admissible region
described by \eqref{eqiq-adm-region}.
By Corollary~\ref{corsuff-opt-tail-expo},
it is clear that the tail exponent in input-queued switches is optimal
under our policy. Moreover, input-queued switches
are not the only network model that satisfies the condition
stated in Corollary~\ref{corsuff-opt-tail-expo}.
For example, consider the independent-set model of a wireless network.
When the underlying interference graph is bipartite,
it is easy to see that the admissible region is characterized
by inequalities of the form $x_i + x_j \leq1$ over
all edges $(i, j)$ of the graph,
and $x_i \leq1$ for isolated nodes $i$.
More generally, when the underlying graph is perfect,
inequality constraints characterizing the
admissible region take the form $\sum_i x_i \leq1$,
where the summation is over all vertices of a clique.
This latter fact follows from a proof of the weak perfect
graph theorem, see, for example, Theorem 12.1.2 in \cite{Matousek2002}.
Thus the incidence matrix $\bR$ has all entries in $\{0, 1\}$,
and the tail exponent under our policy is optimal
for this model.

\subsubsection*{Optimality in input-queued switches}\label{ssseciq}
Here we argue the optimality of our policy for input-queued
switches. As discussed above, the scaling of tail exponent is optimal
under our policy %for any switched networks, and hence
for input-queued switches. We would argue the optimal scaling of the average
total queue size under our policy for input-queued switches. To that end,
as argued in Shah, Tsitsiklis and Zhong \cite{STZopen}, when all
input and output ports
approach critical load,
the average total queue size under
any policy for input-queued switch must scale at least as fast
as $\sqrt{N}/(1-\rho)$,
for any $n$-port switch with $N = n^2$ queues.
For completeness, we include the proof for this lower bound here.
As in Section~\ref{sseciqintro}, we use double indexing.
%
%le5.6 #&#
\begin{lemma}\label{lemqueuesize}
Consider an $n$-port input-queued switch,
with an arrival rate vector $\bolds{\lambda}$.
Suppose that the loads on all input and output ports are $\rho$,
that is, $\sum_{k=1}^n \lambda_{k,\ell} = \sum_{m} \lambda_{\ell,m} = \rho$,
for all $\ell\in\{1, 2, \ldots, n\}$,
where $\rho\in(0,1)$. Consider any policy under which
the queue-size process has a well-defined limiting stationary distribution,
and let this distribution be denoted by ${\bolds{\pi}}_0$.
Then under ${\bolds{\pi}}_0$, we must have
\[
\bbE_{{\bolds{\pi}}_0} \Biggl[\sum_{k,\ell=1}^n
Q_{k,\ell} \Biggr] \geq\frac{n\rho}{2(1-\rho)}.
\]
\end{lemma}
\begin{pf}
We consider the sums of queue sizes at each output port,
that is, the quantities $\sum_{k=1}^n Q_{k, \ell}$ for each $\ell\in
\{1, 2, \ldots, n\}$.
Since at most one packet can depart at each time slot,
$\sum_{k=1}^n Q_{k, \ell}$ stochastically dominates the queue size in
an $M/D/1$ system,
with arrival rate $\rho$ and deterministic service rate $1$.
Therefore, for each $\ell\in\{1, 2, \ldots, n\}$,
\[
\bbE_{{\bolds{\pi}}_0} \Biggl[\sum_{k=1}^n
Q_{k, \ell} \Biggr] \geq\frac
{\rho}{2(1-\rho)}.
\]
Here, $\frac{\rho}{2(1-\rho)}$ is the\vspace*{1pt} expected queue size in steady
state in an $M/D/1$ system.
Summing over $\ell$ gives us the desired bound.
\end{pf}

The optimality in terms of the average total queue size is
a direct consequence of Theorem~\ref{thmmain} and Lemma~\ref{lemqueuesize}.

%
%co5.7 #&#
\begin{corollary}\label{corqueuesize}
Consider the same setup as in Lemma~\ref{lemqueuesize}.
Then in the heavy-traffic limit $\rho\to1$,
our policy is $2$-optimal in terms of the average total queue size.
More precisely, consider the expected total queue size in
the diffusion scale in steady state, that is, $(1-\rho)\widebar{Q}$.
Then
\[
\limsup_{\rho\to1} (1-\rho)\widebar{Q} \leq n
\]
under our policy, and
\[
\liminf_{\rho\to1} (1-\rho)\widebar{Q} \geq\frac{n}{2}
\]
under any other policy.
\end{corollary}
\begin{pf}
Lemma~\ref{lemqueuesize} implies that
\[
\liminf_{\rho\to1} (1-\rho)\widebar{Q} \geq\frac{n}{2}
\]
under any policy. For the upper bound,
note that by Theorem~\ref{thmmain}, under our policy,
\[
\widebar{Q} \leq\frac{J}{2(1-\rho)} + (N+2) K.
\]
For input-queued switches, $J \leq2n$, as remarked in Section~\ref{sseciq},
$N = n^2$ and \mbox{$K = n$}.
Therefore, we have
that under our policy, the expected total queue size satisfies
%
%e33 #&#
\begin{eqnarray}
\widebar{Q} & \leq&\frac{n}{1-\rho} + \bigl(n^2+2\bigr)n.
%& \leq\frac{3 \sqrt{N}}{1-\rho},
\end{eqnarray}
Now consider the steady-state heavy-traffic scaling $(1-\rho)\bQ$.
We have that
%
%e34 #&#
\begin{eqnarray}
(1-\rho)\widebar{Q} & \leq& n + (1-\rho) \bigl(n^2 + 2\bigr)n.
%& \leq\frac{3 \sqrt{N}}{1-\rho},
\end{eqnarray}
The term $(1-\rho) (n^2+2) n$ goes to zero as $\rho\rightarrow1$,
and hence under our policy,
\[
\limsup_{\rho\rightarrow1} (1-\rho)\widebar{Q} \leq n.
\]\upqed
\end{pf}
Our policy is not optimal in terms of the average total queue size,
in general switched networks. In cases where {$J\gg N$},
the moment bounds for the maximum-weight policy gives tighter upper bounds.
For further discussion, see Section~\ref{secconc}.

%%%\iffalse
%%%Theorem~\ref{thmmain}
%%%shows that the expected total queue size under our policy is
%%%bounded by $\frac{J}{2(1-\rho)} + K(N+2)$. For input-queued switches,
%%%$|J| = 2n = 2\sqrt{N}$ and $K = n = \sqrt{N}$. Therefore, we have
%%%that under our policy, the expected total queue-size scales as
%%%{
%%%%
%%%\begin{eqnarray}
%%%\widebar{Q} & \leq\frac{\sqrt{N}}{1-\rho} + (N+2) \sqrt{N}. %
%%%%\nonumber\\
%%%%& \leq\frac{3 \sqrt{N}}{1-\rho},
%%%\end{eqnarray}
%%%%
%%%%for $1/(1-\rho) \geq N+2$.
%%%Now consider the steady-state heavy-traffic scaling $(1-\rho)\bQ$.
%%%We have that
%%%%
%%%\begin{eqnarray}
%%%(1-\rho)\widebar{Q} & \leq\sqrt{N} + (1-\rho) (N+2) \sqrt{N}. %
%%%%\nonumber\\
%%%%& \leq\frac{3 \sqrt{N}}{1-\rho},
%%%\end{eqnarray}
%%%%
%%%The term $(1-\rho) (N+2) \sqrt{N}$ goes to zero as $\rho
%%%and hence under our policy,
%%%%
%%%\[
%%%\limsup_{\rho\rightarrow1} (1-\rho)\widebar{Q} \leq\sqrt{N}.
%%%\]
%%%%
%%%On the other hand, by Lemma~\ref{lemqueuesize}, under any policy,
%%%%
%%%\[
%%%\liminf_{\rho\rightarrow1} (1-\rho)\widebar{Q} \geq\frac{1}{2}
%%%\]
%%%%
%%%Thus, in the heavy-traffic limit, our policy is \emph{optimal}
%%%(up to a constant factor of $2$)
%%%in term of the expected total queue size. }
%%%%Thus, in regime where $1/(1-\rho)$ with
%%%%$\rho= \rho(\blambda)$ scales at a rate faster than $N$, our policy
%%%%results in optimal average queue-size scaling (with universal
%constant
%%%%$3$).
%%%%As a consequence, our policy is within constant factor $3$ of
%optimal
%%%%policy in the so-called heavy traffic scaling.
%%%\fi

%s5.3 #&#
\subsection{Proof of Theorem \texorpdfstring{\protect\ref{thmmain}}{5.2}}\label{ssecproof}

The proof is divided into three parts. The first part describes a
sample-path-wise
relation between $\bQ(\cdot)$ and $\bW(\cdot)$, the residual
workload vector in $\bBN$, which states
that $\bQ(\cdot)$ and $\bW(\cdot)$ differ only by at most a
constant at all
times. Note that this domination is a distribution-free statement. The second
part utilizes this fact to establish the positive recurrence of the
$\bSN$
Markov chain. The third part, as a consequence of the first two
parts, and using Theorem~\ref{thmSFA1}, establishes the quantitative
claims in Theorem~\ref{thmmain}.

\subsubsection*{Part 1. Dominance}
We start by establishing that
the queue sizes $\bQ(\cdot)$ of $\bSN$ are effectively dominated
by the workloads $\bW(\cdot)$ of $\bBN$ at all times.
We state this result formally in Proposition~\ref{propdom}, which is
a consequence
of Lemmas~\ref{lemdom1} and~\ref{lemdom} below.
{%
%le5.8 #&#
\begin{lemma}\label{lemdom1}
Consider the evolution of queue sizes in both $\bBN$ and $\bSN$
networks fed
by identical arrival process. Initially, $\bQ(0) = \bM(0) = \bZero$. Let
$\bW(\tau) = (W_i(\tau))$ denote the amount of unfinished work in
all $N$ queues under the $\bBN$ network at time $\tau$. Then for any
$\tau\geq0$ and $1\leq i\leq N$,
%
%e35 #&#
\begin{equation}
\label{eqd6} %W_i(\tau) & \leq
W_i(\tau) \leq Q_i(\tau)
\leq W_i(\tau) + D_i(\tau), %~\leq M_i(
\end{equation}
where $\bD(\tau) = \Lambda^{\mathrm{SFA}}(\tau) - \bB(\tau)$ is as
described in
Section~\ref{ssecpolicy}.
%where $D_i(\tau) \geq0$ and $\sum_i D_i(\tau) \leq(N+2) K$.
\end{lemma}
\begin{pf}
%To establish the Lemma, we only need to prove \eqref{eqd6} due to
Consider any $i \in\{1, 2, \ldots, N\}$ and $\tau\geq0$.
From \eqref{eqdiscretequeuesinglehop}, in $\bSN$,
%
%e36 #&#
\begin{eqnarray}
\label{eqlemd0} Q_i(\tau) & =& A_i(\tau) -
B_i(\tau) + Z_i(\tau),
\end{eqnarray}
where $Z_i(\tau)$ is the cumulative amount of idling at the $i$th
queue in $\bSN$.
Similarly in $\bBN$,
%
%e37 #&#
\begin{eqnarray}
\label{eqlemd1} W_i(\tau) & =& A_i(\tau) -
\Lambda^{\mathrm{SFA}}_i(\tau) + \hZ_i(\tau),
\end{eqnarray}
where $\hZ_i(\tau)$ is the cumulative amount of idling for the $i$th
queue in $\bBN$.
Since by construction, $\bD(\tau) = \Lambda^{\mathrm{SFA}}(\tau) - \bB
(\tau)$, and
$\bD(\tau) \geq\mathbf{0}$,
we have that
%
%e38 #&#
\begin{eqnarray}
\label{eqlemd2} B_i(\tau) & \leq&\Lambda^{\mathrm{SFA}}_i(
\tau) \leq B_i(\tau) + D_i(\tau).
\end{eqnarray}
By\vspace*{1pt} definition, the instantaneous rate allocation to the $i$th queue satisfies
$\frac{d}{dt^+}\Lambda^{\mathrm{SFA}}_i(t) = 0$ if $W_i(t) = 0$
[equivalently, if
$M_i(t) = 0$]
for any $t \geq0$. Therefore, $\hZ_i(\tau) = 0$, and
$W_i(\tau) = A_i(\tau) - \Lambda^{\mathrm{SFA}}_i(\tau)$. On the other hand,
by Skorohod's map,
%
%e39 #&#
\begin{eqnarray}
\label{eqlemd3} Z_i(\tau) & =& \sup_{0\leq s\leq\tau}
\bigl[B_i(s) - A_i(s) \bigr]^+
\nonumber
\\
& \leq&\sup_{0\leq s\leq\tau} \bigl[\Lambda^{\mathrm{SFA}}_i(s)
- A_i(s) \bigr]^+
\\
& = &\hZ_i(\tau) = 0,\nonumber
\end{eqnarray}
hence for all $i$ and $\tau$, $Z_i(\tau) = 0$, and $Q_i(\tau) =
A_i(\tau) - S_i(\tau)$.
It then follows that
%From \eqref{eqlemd2} and \eqref{eqlemd3}, it follows that
%
%e40 #&#
\begin{eqnarray}
\label{eqlemd4} Q_i(\tau) & =& A_i(\tau) -
B_i(\tau)
\nonumber
\\
& \leq& A_i(\tau) - \Lambda^{\mathrm{SFA}}_i(\tau) +
D_i(\tau)
\\
&=& W_i(\tau) + D_i(\tau)\nonumber
\end{eqnarray}
and
%
%e41 #&#
\begin{eqnarray}
\label{eqlemd5} W_i(\tau) & =& A_i(\tau) -
\Lambda^{\mathrm{SFA}}_i(\tau)
\nonumber\\[-8pt]\\[-8pt]
& \leq& A_i(\tau) - B_i(\tau) = Q_i(\tau).\nonumber
\end{eqnarray}
%
%W_i(\tau) & = A_i(\tau) - \Lambda^{\mathrm{SFA}}_i(\tau) \nonumber\\
%& \leq A_i(\tau) - B_i(\tau) \nonumber\\
%& = Q_i(\tau) ~\nonumber\\
%& \leq A_i(\tau) - \Lambda^{\mathrm{SFA}}_i(\tau) + D_i(\tau) \nonumber
%& = W_i(\tau) + D_i(\tau).
\old{Since the workload at the $i$th queue equals the total amount of
unfinished work for
all of the $M_i(\tau)$ packets waiting at the $i$th queue, and since
each packet
has at most a unit amount of unfinished work, $W_i(\tau) \leq M_i(\tau)$.}
%The proof of Lemma~\ref{lemdom1} is completed by Lemma~\ref{lemdom}
%stated
%below.
Inequalities \eqref{eqlemd4} and \eqref{eqlemd5} together imply
\eqref{eqd6}.
\end{pf}
}
%We state and prove the Lemma~\ref{lemdom} stated below which
%claims that for all $\tau\geq0$, $\bD(\tau)$ is such that $\rho(\bD(
%Formally, we state the following.
%
%le5.9 #&#
\begin{lemma}\label{lemdom}
Let $\bD(\tau)$ be the same as in Lemma~\ref{lemdom1}.
Then, for all $\tau\geq0$, $\rho(\bD(\tau)) \leq N + 2$. In particular,
%
%e42 #&#
\begin{eqnarray}
\label{eqDsum} \sum_i D_i(\tau) &
\leq& K(N+2)\qquad\mbox{where }K = \max_{\bolds{\sigma}
\in\sS} \sum
_i \sigma_i.
\end{eqnarray}
\end{lemma}
\begin{pf}
This result is established as follows. First, observe that $\bD(0) =
\bZero$
and therefore $\rho(\bD(0)) = 0$. Next, we show that
$\rho(\bD(\tau+1)) \leq\rho(\bD(\tau)) + 1$. That is, $\rho(\bD
(\cdot))$ can
at most increase by $1$ in each time slot. And finally, we show that~$\rho(\bD(\cdot))$ cannot
increase once it exceeds $N+1$. That is, if $\rho(\bD(\tau)) \geq
N+1$, then
$\rho(\bD(\tau+1))\leq\rho(\bD(\tau))$. This will complete the proof.

We start by establishing that $\rho(\bD(\cdot))$ increases
by at most $1$ in unit time. By definition,
\begin{eqnarray}\label{eqd1}
\bD(\tau+ 1) &=& \Lambda^{\mathrm{SFA}}(\tau+1) - \bB(\tau+1)\nonumber
\\
& =& \Lambda^{\mathrm{SFA}}(\tau) - \bB(\tau) + \bigl(\Lambda^{\mathrm{SFA}}(
\tau+1)-\Lambda^{\mathrm{SFA}}(\tau) - d\bB(\tau) \bigr)\nonumber
\nonumber\\[-8pt]\\[-8pt]
& =& \bD(\tau) + d\Lambda^{\mathrm{SFA}}(\tau) - d\bB(\tau)\nonumber
\\
& =& \bigl(\bD(\tau) - d\bB(\tau) \bigr) +d\Lambda^{\mathrm{SFA}}(\tau),\nonumber
\end{eqnarray}
where $d\Lambda^{\mathrm{SFA}}(\tau) = \Lambda^{\mathrm{SFA}}(\tau+1)
-\Lambda^{\mathrm{SFA}}(\tau)$.
As remarked earlier, $d\bB(\tau) \leq\bD(\tau)$ componentwise.
Therefore, by \eqref{eqrhosubadd} it follows that
\[
\rho\bigl(\bD(\tau+1) \bigr) \leq\rho\bigl(\bD(\tau)-d\bB(\tau)
\bigr) +
\rho\bigl(d\Lambda^{\mathrm{SFA}}(\tau)\bigr).
\]
Note that $\rho(d\Lambda^{\mathrm{SFA}}(\tau)) \leq1$ because the instantaneous
service rate under SFA is always admissible.
Since $\bD(\tau) \geq\bD(\tau) - d\bB(\tau) \geq\mathbf{0}$,
any feasible solution to $\PRIMAL(\bD(\tau) )$
is also feasible to $\PRIMAL(\bD(\tau) - d\bB(\tau) )$,
and hence
\[
\rho\bigl(\bD(\tau)-d\bB(\tau) \bigr) \leq\rho\bigl(\bD(\tau) \bigr).
\]
It follows that
%
%e43 #&#
\begin{eqnarray}
\label{eqd2} \rho\bigl(\bD(\tau+1) \bigr) & \leq&\rho\bigl(\bD(\tau
)\bigr) + 1.
\end{eqnarray}
Next, we shall argue that if $\rho(\bD(\tau)) \geq N+1$, then
$\rho(\bD(\tau+1))\leq\rho(\bD(\tau))$. To that end,\vspace*{1pt} suppose that
$\rho(\bD(\tau)) \geq N+1$. Now $\frac{1}{\rho(\bD(\tau))} \bD
(\tau) \in\hull{\sS}$.
Note that $\hull{\sS}$ is a convex set in a $N$-dimensional space
with extreme points
contained in $\sS$. Therefore, by Carath\'{e}odory's theorem, $\frac
{1}{\rho(\bD(\tau))} \bD(\tau)$
can be\vspace*{2pt} written as a convex
combination of at most $N + 1$ elements in $\sS$. That is,
there exists $\alpha_k \geq0$ with $\sum_{k=1}^{N+1} \alpha_k = 1$,
and $\bolds{\sigma}^k \in\sS$,
$k \in\{1, 2, \ldots, N+1\}$, such that
%
%e44 #&#
\begin{eqnarray}
\label{eqd3} \frac{1}{\rho(\bD(\tau))} \bD(\tau) & =& \sum
_{k=1}^{N+1} \alpha_{k} \bolds{
\sigma}^k.%,~~\bolds{\sigma}^k \in\sS, ~\alpha_k \geq0, ~~
\end{eqnarray}
Therefore, there exists some $k^* \in\{1, 2, \ldots, N+1\}$,
such that $\alpha_{k^*} \geq1/\break (N+1)$. Since $\rho(\bD(\tau)) \geq N+1$,
$\rho(\bD(\tau))\alpha_{k^*} \geq1$. That is, $\bD(\tau)$ can be
written as a nonnegative combination of
elements from $\sS$ with one of them, $\bolds{\sigma}^{k^*}$, having an
associated coefficient
that satisfies $\rho(\bD(\tau)) \alpha_{k^*} \geq1$,
as required. In this case, we have
%
%e45 #&#
\begin{equation}
\bD(\tau) - \bolds{\sigma}^{k^*} = \sum_{k=1, k\neq k^*}^{N+1}
{\rho\bigl(\bD(\tau)\bigr)}\alpha_k \bolds{\sigma}^k +
\bigl({\rho\bigl(\bD(\tau)\bigr)}\alpha_{k^*}-1\bigr)\bolds{
\sigma}^{k^*}.
\end{equation}
Therefore,
%
%e46 #&#
\begin{equation}
\label{eqd4} \rho\bigl(\bD(\tau)-\bolds{\sigma}^{k^*} \bigr) \leq\rho
\bigl(\bD(\tau)\bigr) - 1.
\end{equation}
Our scheduling policy chooses such a schedule, that is, $\bolds{\sigma}^{k^*}$;
that is, $d\bB(\tau) = \bolds{\sigma}^{k^*}$. Therefore,
%
%e47 #&#
\begin{eqnarray}
\bD(\tau+1) & =& \bD(\tau) - \bolds{\sigma}^{k^*}+ d
\Lambda^{\mathrm{SFA}}(\tau).
\end{eqnarray}
By another application of \eqref{eqrhosubadd} it follows that
%
%e48 #&#
\begin{eqnarray}
\rho\bigl(\bD(\tau+1)\bigr) & \leq&\rho\bigl(\bD(\tau) - \bolds{
\sigma}^{k^*} \bigr) + \rho\bigl(d\Lambda^{\mathrm{SFA}}(\tau)\bigr)
\nonumber
\\
& \leq&\rho\bigl(\bD(\tau) \bigr) - 1 + 1,
\\
& =& \rho\bigl(\bD(\tau) \bigr),\nonumber
\end{eqnarray}
where again we have used the fact that $\rho(d\Lambda^{\mathrm
{SFA}}(\tau))\leq1$,
due to the feasibility of SFA
policy and \eqref{eqd4}. This establishes that $\rho(\bD(\tau))
\leq N+2$ for all $\tau\geq0$.
That is, for each $\tau\geq0$, there exists $\alpha_{\bolds{\sigma}} \geq0$
for all $\bolds{\sigma}\in\sS$,
$\sum_{\bolds{\sigma}} {\rho(\bD(\tau))}\alpha_{\bolds{\sigma}} \leq
N+2$ and
%
%e49 #&#
\begin{eqnarray}
\bD(\tau) & \leq&\sum_{\bolds{\sigma}} \alpha_{\bolds{\sigma}}
\bolds{\sigma}.
\end{eqnarray}
Therefore,
%
%e50 #&#
\begin{eqnarray}
\label{eqd5} \sum_i D_i(\tau) & =&
\bD(\tau) \bdot\bOne\leq\sum_{\bolds{\sigma}} {\rho\bigl(\bD(
\tau)\bigr)}\alpha_{\bolds{\sigma}} \bolds{\sigma}\bdot\bOne
\nonumber\\[-8pt]\\[-8pt]
& \leq&\biggl(\sum_{\bolds{\sigma}} {\rho\bigl(\bD(\tau)\bigr)}
\alpha_{\bolds{\sigma}} \biggr) \biggl(\max_{\bolds{\sigma}\in\sS} \sum
_i \sigma_i \biggr)
\leq (N+2) K,\nonumber
\end{eqnarray}
where $K = \max_{\bolds{\sigma}\in\sS} \sum_i \sigma_i$. This completes
the proof of
Lemma~\ref{lemdom}.
\end{pf}

Lemmas~\ref{lemdom1} and~\ref{lemdom} together imply the following proposition.
%
%pr5.10 #&#
\begin{proposition}\label{propdom}
Let $\bQ(\cdot)$, $\bW(\cdot)$ and $\bM(\cdot)$ be as in Lemma
\ref{lemdom1}. Then
%
%e51 #&#
\begin{equation}
\label{eqdomfinal} \sum_{i=1}^N
Q_i(\tau) \leq\sum_{i=1}^N
W_i(\tau) + K(N+2) \leq\sum_{i=1}^N
M_i(\tau) + K(N+2),
\end{equation}
where $K = \max_{\bolds{\sigma}\in\sS} (\sum_{i=1}^N \sigma_i )$.
\end{proposition}
\begin{pf}
We obtain the bounds of \eqref{eqdomfinal} by summing inequality
\eqref{eqd6} over $i \in\{1, 2, \ldots, N\}$,
and using bound \eqref{eqDsum}.
\end{pf}

\subsubsection*{Part 2. Positive recurrence}
We start by defining the Markov chain describing the system evolution
under the policy of interest. There are essentially two systems
that evolve in a coupled manner under our policy: the virtual
bandwidth-sharing network $\bBN$ and the switched network
$\bSN$ of interest. These two networks are fed by the
same arrival processes which are exogenous and Poisson (and
hence Markov). The virtual system $\bBN$ has a Markovian\vadjust{\goodbreak}
state consisting of the packets whose services are not completed,
represented by the vector $\bM(\cdot)$, and their residual services.
The residual services of $M_i(\cdot)$ packets queued on
route $i$ can be represented by a nonnegative, finite
measure $\mu_i(\cdot)$ on $[0,1]$: unit mass is placed
at each of the points $0\leq s_1,\ldots, s_{M_i(t)} \leq1$
if the unfinished work of $M_i(t)$ packets are given by
$0 < s_1,\ldots, s_{M_i(t)} \leq1$.

We now consider a Markovian description of the
network $\bSN$ in discrete time: let $\bX(\tau)$ be the state of the
system defined as
%
%e52 #&#
\begin{equation}
\label{eqsnmarkov} \bX(\tau) = \bigl(\bM(\tau), \bmu(\tau), \bQ(\tau),
\bD(\tau)
\bigr),
\end{equation}
where $(\bM(\tau), \bmu(\tau))$ represents the state of $\bBN$ at
time $\tau$,
$\bQ(\tau)$ is the vector of queue sizes in $\bSN$ at time $\tau$
and $\bD(\tau)$
is the ``difference'' vector maintained by the scheduling policy for
$\bSN$,
as described in Section~\ref{ssecpolicy}. Observe that the state \mbox{$\bX
(\tau+1)$} is a function of the previous state $\bX(\tau)$ and the
independent random arrival times occurring in the time interval $(\tau,\tau+1)$ according to our Poisson process. This ensures conditional
independence between $\bX(\tau+1)$ and $\bX(\tau-1)$ given $\bX
(\tau)$. So, by standard arguments,
$\bX(\cdot)$ is Markov and, indeed, strong Markov.

We now define the state space $\sX$ of the Markov chain $\bX(\cdot)$:
%%%\iffalse
%%%Informally speaking, $\sX$ will consist of points
%%%that not only can be reached from the zero state $\bZero$,
%%%but also can reach $\bZero$, in finite time.
%%%More precisely, first note that $\sX$ is a subset of the product
%space
%%%\fi
%
\[
{\sX=\IntegersP^N \times\cM\bigl([0,1]\bigr)^N \times
\IntegersP^N \times\mathcal{D}},
\]
where $\cM([0,1])$ is the space of all nonnegative, finite measures
on $[0,1]$
{and where $\mathcal{D}=(N+2)\cdot\mathcal{C}$ is the admissible
region $\mathcal{C}$ expanded by a multiplicative factor $N+2$.
The set $\mathcal{D}$ is exactly the set of vectors $\mathbf{d} \in\bR_+^N$ for which $\rho(\mathbf{d})\leq N+2$; cf. \eqref
{eqdef-load-2}. Thus, by Lemma
\ref{lemdom}, our process $\bD(\cdot)$ can never leave the set
$\mathcal{D}$. }

We endow $\cM([0,1])$ with the weak topology, which is induced by the
Prohorov's metric.
This results in a complete and separable metric (Polish) space.
The set $\mathcal{D}$ is a closed convex subset of $\bR_+^N$. We
endow $\IntegersP$ and $\mathcal{D}$ with the obvious metrics (e.g.,
$\ell_1$).
The entire product space is endowed with the metric that is the maximum
of metrics
on component spaces. The resulting product space is Polish,
on which a Borel $\sigma$-algebra, $\mathcal{B}_\sX$, can be defined.
%%%\iffalse
%%%Then, $\sX$ consists of points $\mathbf{x}$ in this product space
%%%such that
%%%%
%%%\[
%%%\Prob_{\mathbf{x}}\left(T_{\bZero} < \infty\right) = 1, \quad
%%%\mbox{and } \quad
%%%\Prob_{\bZero}\left(T_{\mathbf{x}} < \infty\right) = 1.
%%%\]
%%%%
%%%Here $\Prob_{\mathbf{x}}(A) \equiv\Prob(A \mid\bX(0) = \mathbf{x})$
%%%denotes probability conditioned on the initial condition $\bX(0) =
%%%\mathbf{x}$,
%%%and for a measurable set $A$, and $T_{A} = \inf\{\tau\geq1 \dvtx  \bX
%%%(\tau) \in A\}$
%%%denotes the return time to $A$.
%%%\fi

We remark that the Markov chain $\bX(\cdot)$ need not be recurrent
(nor neighborhood recurrent) for all states in $\sX$.
However, we can start our Markov chain from any state $x\in\sX$ and
it will hit state $\bZero$ in finite expected time. We can then prove
that our Markov chain is positive Harris recurrent. The resulting
stationary measure defines the subset of $\sX$ for which $\bX(\cdot
)$ is recurrent.

Given the Markovian description $\bX(\tau)$ of $\bSN$,
we establish its positive Harris recurrence in the following lemma.

%
%le5.11 #&#
\begin{lemma}\label{LEMPOSREC}
Consider a switched network $\bSN$ with a strictly admissible arrival
rate vector $\bolds{\lambda}$,
with $\rho(\bolds{\lambda}) < 1$. Suppose that at time $0$ the
system is empty.
Let $\bX(\cdot)$ be as defined in equation \eqref{eqsnmarkov}.
Then $\bX(\cdot)$ is positive Harris recurrent {and ergodic}.\vadjust{\goodbreak}
\end{lemma}
The proof of the lemma is technical, and is deferred to Appendix~\ref{apdxposrec}.
The idea is that the evolution of $\bBN$ is not affected by $\bSN$,
and that $\bBN$ is, on its own, positive recurrent.
Hence, starting from any initial state, the Markov process $(\bM(\cdot
), \bmu(\cdot))$
that describes the evolution of $\bBN$, reaches the null state, that is,
$(\bM(\cdot), \bmu(\cdot)) = \bZero$ at some finite expected time.
Once $\bBN$ reaches the null state, it stays at this state for an
arbitrarily large amount of time with positive probability.
By our policy, $\bQ(\cdot)$ and $\bD(\cdot)$ can be driven to
$\bZero$
within this time interval. This establishes that $\bX(\cdot)$ reaches
the null state in finite expected time and that $\bX(\cdot)$
is positive recurrent.

\subsubsection*{Part 3. Completing the proof}
The positive recurrence of the Markov chain $\bX(\cdot)$
implies that it possesses a unique stationary distribution. %and that
%is ergodic.
Let $\widebar{W} = \bbE_{{\bolds{\pi}}} [\sum_{i=1}^N W_i ]$,
where, similarly to Proposition~\ref{propworkload}, $W_i$ is the
steady-state workload on queue $i$
in $\bBN$. %Define $\widebar{M}$ similarly.
By ergodicity, the time average of the total queue size
equals the expected total queue size in steady state, that is, $\widebar{Q}$,
and similarly for~$\widebar{W}$. Therefore, by Proposition~\ref{propdom},
\[
\widebar{Q} \leq\widebar{W} + K(N+2).
\]
By Proposition~\ref{propworkload},
\[
\widebar{W} = \frac{1}{2} \Biggl(\sum_{j=1}^J
\frac{\tilde{\rho}_j}{1-\tilde{\rho}_j} \Biggr).
\]
Thus
\[
\widebar{Q} \leq\widebar{W} + K(N+2) = \frac{1}{2} \Biggl(\sum
_{j=1}^J \frac{\tilde{\rho}_j}{1-\tilde{\rho}_j} \Biggr) + K(N+2).
\]

We now establish the tail exponent in \eqref{eqthmmain2}.
By Proposition~\ref{propdom},
\[
\sum_{i=1}^N W_i(\tau) \leq
\sum_{i=1}^N Q_i(\tau) \leq
\sum_{i=1}^N W_i(\tau) +
K(N+2),
\]
deterministically and for all times $\tau$.
Since $K(N+2)$ is a constant, $\sum_{i=1}^N Q_i(\cdot)$ and
$\sum_{i=1}^N W_i(\cdot)$ have the same tail exponent in steady state.
By Proposition~\ref{propworkload},
the tail exponent $\beta(\bW)$ of $\sum_{i=1}^N W_i$ in steady state
is given by $-\theta^*$,
where $\theta^*$ is the unique positive solution of the equation
$\rho(e^{\theta} - 1) = \theta$,
so
\[
\beta(\bQ) = \beta(\bW) = -\theta^*.
\]
\section{Discussion}\label{secconc}

We present a novel scheduling policy for a general single-hop switched
network model. The
policy, in effect, emulates the so-called Store-and-forward (SFA)
continuous-time
bandwidth-sharing policy. The insensitivity property of SFA along with
the relation of its stationary distribution with that of a multi-class queuing
network leads to the explicit characterization of the stationary\vadjust{\goodbreak}
distribution of
queue sizes induced by our policy. This allows us to establish the optimality
of our policy in terms of tail exponent for a large class of switched networks,
including input-queued switches, and the independent-set model of
wireless networks
when the underlying interference graph is perfect, and
that with \mbox{respect} to the average total queue size for a class of switched
networks, including the input-queued switches. As a consequence, this settles
a conjecture stated in~\cite{STZopen}.
On the technical end, a key contribution of the paper is creating a
discrete-time
scheduling policy from a continuous-time rate allocation policy,
and this may be of independent interest in other domains of applications.
We also remark that the idea of designing a discrete-time policy by emulating
a continuous-time policy is not new; for example, similar emulation schemes
have appeared in \cite{DKS90,GP93}. Our emulation scheme is novel
in that it captures the switched network structure where queues may be
served simultaneously. This simultaneity of service is absent from
earlier models.

The switched network model considered here requires the arrival processes
to be Poisson. However, this is not a major restriction, due to a
\textit{Poissonization} trick considered, for example, in \cite{EMPS} and
\cite{JS}: all arriving packets are first passed through a ``regularizer,''
which emits out packets according to a Poisson process with a rate
that lies between the arrival rate and the network capacity. This
leads to the arrivals being effectively Poisson, as seen by the system
%%seen by the system look like a Poisson process
with a somewhat higher rate---by choosing the rate of ``regularizer'' so that
the effective gap to the capacity, that is, $(1-\rho)$, is decreased
by factor $2$.

The scheduling policy that we propose is not optimal for general
switched networks.
For example, in the independent-set model of ad-hoc wireless networks,
there are as many constraints as
the number of edges in the interference graph, which is often
much larger than the number of nodes.
Under our policy, the average total queue size would scale with
the number of edges, whereas maximum-weight policy achieves a scaling
with the number of nodes.

There are many possible directions for future research.
One direction is the search for low-complexity and optimal scheduling policies.
In the context of input-queued switches,
our policy has a complexity that is exponential
in $N$, the number of queues, because one has to compute the sum of
exponentially many terms
at every time instance. This begs the question
of finding an optimal policy with polynomial complexity in $N$.
One candidate is the MW-$\alpha$ policy, $\alpha> 0$, which has
polynomial complexity,
but its optimality appears difficult to analyze.
Another possible candidate could be, as discussed in the \hyperref[secintro]{Introduction},
a randomized version of proportional fairness.
The relationship between SFA and
proportional fairness is explored in \cite{walton},
where it was formally established that SFA converges to proportional fairness
under the heavy-traffic limit, in an appropriate sense.
The question remains whether (a version of) proportional fairness is optimal
for input-queued switches.

Another interesting direction to pursue has to do
the analysis of different limiting regimes.
We are interested in two limits: $N \to\infty$,
and $\rho\to1$, where $N$ is the number of queues,
and $\rho$ is the system load. Again, take the example
of input-queued switches. In this paper,
we have considered the heavy-traffic limit, that is, $\rho\to1$,
and show that our policy is optimal.
However, if we take the limit $N \to\infty$, while keeping $\rho$
fixed, then the average total queue-size scales as $N^{3/2}$,
whereas maximum-weight policy produces a bound of $N$.
A more interesting question is in the regime where $(1-\rho)\sqrt{N}$
remain bounded, and where $N \to\infty$. In this regime,
under our policy, under the maximum-weight policy, and under the
batching policy
in \cite{NeelyModiano}, the average total queue sizes
all scale as $O(N^{3/2})$. In contrast, the scaling conjectured in
\cite{STZopen}
is $O(N)$.
It is therefore of interest to see whether the $3/2$ barrier can be broken.
In \cite{STZ25}, the authors
device a policy that achieves $N^{\gamma}$ scaling,
for some $\gamma\in[1, 3/2)$.

%The question of determining optimal policy in terms of average
%queue-sizes for {\em any} switch network remains an outstanding open
%problem. One possibility is to utilize a different insensitive
%bandwidth
%allocation policy for the purpose of emulation that leads to better
%dependence on the underlying structure when $J$ is large.

%sA #&#
\begin{appendix}\label{app}
%sB #&#
\section{Properties of SFA}\label{apdxpsn}

This section {proves} results stated in
Section~\ref{secinsensitive}, specifically Theorem~\ref{thmSFA1},
Propositions~\ref{propSFA2},~\ref{propMC} and~\ref{propworkload}.
First, we note that Propositions~\ref{propSFA2}
and~\ref{propMC} are fairly easy consequences of Theorem~\ref{thmSFA1},
and their proofs are included for completeness.
We then prove Proposition~\ref{propworkload}.
Theorem~\ref{thmSFA1} follows from the work of Zachary \cite
{zachary}.

\begin{pf*}{Proof of Proposition~\ref{propSFA2}}
\old{
First note that to show that $\tilde{{\bolds{\pi}}}$ is a probability measure,
it suffices to establish \eqref{eqsfaeqmc}.
This is because ${\bolds{\pi}}$ is a probability measure,
so once \eqref{eqsfaeqmc} is established,
we can sum over $L$ to get
\[
\sum_{\tilde{\mathbf{m}}\in\IntegersP^{|\cK|}} \tilde{{\bolds{\pi
}}}(\tilde{
\mathbf{m}}) = 1.
\]
}
To verify \eqref{eqsfaeqmc}, we can calculate both sides of the
equation directly.
Note that by definition, $\tilde{m}_j = \sum_{i \dvtx  j \in i} \tilde
{m}_{ji}$, so
%
%eB.1 #&#
\begin{equation}
\label{eqsfaeqmcR} \tilde{{\bolds{\pi}}} \Biggl( \Biggl\{\tilde{\mathbf
{m}}\dvtx  \sum
_{j=1}^J \tilde{m}_j = L
\Biggr\} \Biggr) = \tilde{{\bolds{\pi}}} \biggl( \biggl\{\tilde{\mathbf
{m}}\dvtx  \sum
_{(j,i) \in\cK} \tilde{m}_{ji} = L \biggr\}
\biggr).
\end{equation}
On the other hand,
%
%eB.2 #&#
%eB.3 #&#
%eB.4 #&#
%eB.5 #&#
%eB.6 #&#
\begin{eqnarray}
\qquad&& {\bolds{\pi}} \Biggl( \Biggl\{\mathbf{m}\dvtx  \sum_{i=1}^N
m_i = L \Biggr\} \Biggr)
\nonumber
\\
&&\qquad= \sum_{{\mathbf{m}\in\bbZ_+^{|I|}}} {\bbI\Biggl[\sum
_{i=1}^N m_i = L \Biggr]}
\frac{\Phi(\mathbf{m})}{\Phi} \prod_{i=1}^N
\lambda_i^{m_i} \label{eqsfaeqmcL1}
\\
&&\qquad= \frac{1}{\Phi} \sum_{{\mathbf{m}\in\bbZ_+^{|I|}}} {\bbI
\Biggl[
\sum_{i=1}^N m_i = L \Biggr]}
\nonumber\\[-8pt]\label{eqsfaeqmcL2}\\[-8pt]
&&\quad\qquad{}\times \sum_{\tilde{\mathbf{m}} \in U(\mathbf{m})} \prod_{i=1}^N
\lambda_i^{m_i} \prod_{j=1}^J
\biggl(\pmatrix{\tilde{m}_j \vspace*{2pt}\cr\tilde{m}_{ji}\dvtx  i \ni j}\prod
_{i\dvtx  j\in i} \biggl(\frac{R_{ji}}{C_j} \biggr)^{\tilde{m}_{ji}}
\biggr)\nonumber
\\
&&\qquad= \frac{1}{\Phi} \sum_{{\mathbf{m}\in\bbZ_+^{|I|}}} \sum
_{\tilde{\mathbf{m}} \in U(\mathbf{m})} {\bbI\Biggl[\sum_{i=1}^N
m_i = L \Biggr]}
\nonumber\\[-8pt]\label{eqsfaeqmcL3} \\[-8pt]
&&\quad\qquad{}\times \prod_{j=1}^J
\biggl(\pmatrix{\tilde{m}_j \vspace*{2pt}\cr\tilde{m}_{ji}\dvtx  i \ni j}\prod
_{i\dvtx  j\in i} \biggl(\frac{R_{ji}\lambda_i}{C_j} \biggr)^{\tilde{m}_{ji}}
\biggr) \nonumber
\\
&&\qquad= \frac{1}{\Phi} \sum_{{\tilde{\mathbf{m}}\in\bbZ_+^{|\cK|}}
}{\bbI\biggl[\sum
_{(j,i) \in\cK}\tilde{m}_{ji} = L \biggr]}
\nonumber\\[-8pt]\label{eqsfaeqmcL4} \\[-8pt]
&&\quad\qquad{}\times\prod
_{j=1}^J \biggl(\pmatrix{\tilde{m}_j
\vspace*{2pt}\cr\tilde{m}_{ji}\dvtx  i \ni j}\prod_{i\dvtx  j\in i}
\biggl(\frac{R_{ji}\lambda_i}{C_j} \biggr)^{\tilde{m}_{ji}} \biggr)\nonumber
\\
&&\qquad= \tilde{{\bolds{\pi}}} \biggl( \biggl\{\tilde{\mathbf{m}}\dvtx
\sum
_{(j,i) \in\cK} \tilde{m}_{ji} = L \biggr\} \biggr)
\label{eqsfaeqmcL5}.
\end{eqnarray}
Equality \eqref{eqsfaeqmcL1} follows from the definition of ${\bolds{\pi}}$
given in \eqref{eqsfameasure},
\eqref{eqsfaeqmcL2} follows from the definition of $\Phi(\mathbf{m})$
given in \eqref{dfsfa1},
\eqref{eqsfaeqmcL3} follows from the fact that for $\tilde{\mathbf{m}}
\in
U(\mathbf{m})$,
$\sum_{j\dvtx  j\in i} \tilde{m}_{ji} = m_i$ for all $i \in\cI$,
\eqref{eqsfaeqmcL4} follows from the fact that
\begin{eqnarray*}
&& \sum_{{\tilde{\mathbf{m}} \in\bbZ_+^{|\cK|}, \mathbf{m}\in\bbZ
_+^{|I|}}}\bbI\Biggl[\sum
_{i=1}^N m_i = L, \sum
_{j\dvtx  j\in i} \tilde{m}_{ji} = m_i \Biggr]
=
\bbI\biggl[\sum_{(j,i) \in\cK}\tilde{m}_{ji} = L
\biggr]
\end{eqnarray*}
and \eqref{eqsfaeqmcL5} follows from the definition of $\tilde{{\bolds
{\pi}}}$
given in \eqref{eqmcmeasure}.
Equations \eqref{eqsfaeqmcR} and \eqref{eqsfaeqmcL5} together
establish \eqref{eqsfaeqmc}.
\end{pf*}

\old{
From \eqref{eqsummation}, \eqref{eqcombine1}, and \eqref{eqcombine2},
we see that under $\tilde{{\bolds{\pi}}}$, the marginals $\tilde{m}_j$ are
independently distributed,
with
\[
\tilde{{\bolds{\pi}}} \biggl( \biggl\{\tilde{\mathbf{m}}\dvtx  \sum
_{i\dvtx  j \in i} \tilde{m}_{ji} = \tilde{m}_j
\biggr\} \biggr) = \frac{1}{1 - (\sum_{i\dvtx j\in i} \rho_i)/{C_j}}
\biggl(\frac{\sum_{i\dvtx j\in i} \rho_i}{C_j}
\biggr)^{\tilde{m}_j}.
\]
Thus, on the one hand,
%
%eB.7 #&#
%eB.8 #&#
%eB.9 #&#
%eB.10 #&#
\begin{eqnarray}
\tilde{{\bolds{\pi}}} \biggl( \biggl\{\tilde{\mathbf{m}}\dvtx  \sum
_{j \in J} \tilde{m}_j = L \biggr\} \biggr) & = & \sum
_{\substack{(\tilde{m}_j)_{j\in J}, \\ \sum_{j \in J} \tilde{m}_j =
L}} \prod_{j \in J}
\frac{1}{1 - {\sum_{i\dvtx j\in i} \rho_i}/{C_j}} \biggl(\frac{\sum
_{i\dvtx j\in i} \rho_i}{C_j} \biggr)^{\tilde{m}_j}
\\
& = & \frac{1}{\Phi}\sum_{\substack{(\tilde{m}_j)_{j\in J}, \\ \sum
_{j \in J} \tilde{m}_j = L}} \prod
_{j \in J} \biggl(\frac{\sum_{i\dvtx j\in i} \rho_i}{C_j} \biggr)^{\tilde{m}_j}.
\end{eqnarray}
On the other hand,
}

\begin{pf*}{Proof of Proposition~\ref{propMC}}
We can verify \eqref{eqmcindep} directly.
Indeed,
%
%eB.11 #&#
%eB.12 #&#
%eB.13 #&#
\begin{eqnarray}
\qquad && \tilde{{\bolds{\pi}}} \bigl(\{\tilde{m}_j = L_j \dvtx
j = 1, 2, \ldots, J\} \bigr)
\nonumber
\\
&&\qquad= \frac{1}{\Phi} \sum_{{\tilde{\mathbf{m}}\in\bbZ_+^{|\cK|}}}
{\bbI\Biggl[
\sum_{i=1}^N\tilde{m}_{ji} =
L_j \Biggr]} \prod_{j=1}^J
\biggl(\pmatrix{L_j \vspace*{2pt}\cr\tilde{m}_{ji}\dvtx  i \ni j}\prod
_{i\dvtx  j\in i} \biggl(\frac
{R_{ji}\lambda_i}{C_j} \biggr)^{\tilde{m}_{ji}} \biggr)
\label
{eqsummation}
\\
%= & ~~\sum_{\tilde{m}_j \in\Integers, j \in\cJ} \sum_{\substack{
&&\qquad= \frac{1}{\Phi} \prod_{j=1}^J \biggl(\sum
_{i\dvtx  j \in i}\frac
{R_{ji}\lambda_i}{C_j} \biggr)^{L_j}
\label{eqcombine1}
\\
&&\qquad= \prod_{j=1}^J \biggl(
\frac{C_j - \sum_{i\dvtx i \ni j}R_{ji} \lambda
_i}{C_j} \biggr) \biggl(\sum_{i\dvtx  j \in i}
\frac{R_{ji}\lambda
_i}{C_j} \biggr)^{L_j} \label{eqcombine2}
\\
&&\qquad= \prod_{j=1}^J (1 - \tilde{
\rho}_j )\tilde{\rho}_j^{L_j}.
\nonumber
\end{eqnarray}
Equality\vspace*{1pt} (\ref{eqsummation}) follows from the definition of $\tilde
{{\bolds{\pi}}
}$ in \eqref{eqmcmeasure}.
%for fixed $\tilde{m}_j, j \in\cJ$,
%we first sum over all $\tilde{\mathbf{m}} \in\IntegersP^{|\cK|}$ such
%that
%$\sum_{i\dvtx  j \in i} \tilde{m}_{ji} = \tilde{m}_j$ for all $j \in\cJ$,
%and then we sum over all $\tilde{m}_j \in\Integers, j \in\cJ$.
Equality \eqref{eqcombine1} collects all terms in the Newton expansion
of the term
$ (\sum_{i\dvtx  i \ni j}\frac{R_{ji} \lambda_i}{C_j} )^{L_j}$.
Equality \eqref{eqcombine2} follows from the definition of $\Phi$.
%is also a collection of Newton expansion.
%Equality \eqref{eqtaylor} uses the Taylor expansion of the function
%$f(x) = 1/(1-x)$
%for $x$ with $ |x| < 1$.
\end{pf*}

\begin{pf*}{Proof of Proposition~\ref{propworkload}}
Consider $\sum_{i=1}^N M_i$, the total number of packets waiting in
$\bBN$, in steady state.
By Propositions~\ref{propSFA2} and~\ref{propMC},
$\sum_{i=1}^N M_i$ has the same distribution as the sum of $J$
geometric random variables,
with parameters $1-\tilde{\rho}_1, \ldots, 1-\tilde{\rho}_J$. Hence,
\[
\bbE\Biggl[\sum_{i=1}^N M_i
\Biggr] = \sum_{j=1}^J \frac{\tilde{\rho}_j}{1-\tilde{\rho}_j}.
\]
By Theorem~\ref{thmSFA1}, the individual residual workload in steady
state is independent
from the number of packets in the network, and is uniformly distributed
on [0, 1].
Thus
\[
\bbE\Biggl[\sum_{i=1}^N W_i
\Biggr] = \frac{1}{2}\bbE\Biggl[\sum_{i=1}^N
M_i \Biggr] = \frac{1}{2}\sum_{j=1}^J
\frac{\tilde{\rho}_j}{1-\tilde{\rho}_j}.
\]
This establishes equation \eqref{eqmean-workload}.

To establish equation \eqref{eqtail-workload},
consider the following interpretation of\break  $\sum_{i=1}^N W_i$,
the total residual workload in steady state.
By Theorem~\ref{thmSFA1}, $\sum_{i=1}^N W_i$
has the same distribution as $\sum_{\ell=1}^{M} U_{\ell}$,
where $M = \sum_{i=1}^N M_i$, and $U_{\ell}$
are i.i.d. uniform random variables on $[0, 1]$,
all independent from $M$.
We first establish that
%
%eB.14 #&#
\begin{equation}
\label{eqtail-workload-upper} \limsup_{L \to\infty} \frac{1}{L} \log
\Prob
\Biggl(\sum_{\ell
=1}^M U_{\ell} \geq
L \Biggr) \leq-\theta^*,
\end{equation}
where $\theta^*$ is the unique \emph{positive} solution of the equation
$\rho(e^{\theta} - 1) = \theta$.
By Markov's inequality, for any $\theta> 0$, we have
\begin{eqnarray*}
\Prob\Biggl(\sum_{\ell=1}^M
U_{\ell} \geq L \Biggr) & \leq& \exp(-\theta L) \bbE\Biggl[\exp\Biggl(
\theta\sum_{\ell
=1}^M U_{\ell}
\Biggr) \Biggr]
\\
& = & \exp(-\theta L) \bbE\Biggl[ \bbE\Biggl[ \exp\Biggl(\theta\sum
_{\ell=1}^M U_{\ell} \Biggr)\Bigg| M \Biggr]
\Biggr]
\\
& = & \exp(-\theta L) \bbE\biggl[ \biggl(\frac{e^{\theta}-1}{\theta
}
\biggr)^M \biggr].
\end{eqnarray*}
For notational convenience, let $x = \frac{e^{\theta}-1}{\theta}$.
We now consider the term $\bbE[x^M]$.
Let $\widetilde{M}_j$ be independent geometric random variables
with parameter $1-\tilde{\rho}_j$, $j = 1, 2, \ldots, J$.
Then $M$ is distributed as $\sum_{j=1}^J \widetilde{M}_j$.
Thus
\[
\bbE\bigl[x^M\bigr] = \bbE\bigl[x^{\sum_{j=1}^J \widetilde{M}_j} \bigr
] = \prod
_{j=1}^J \bbE\bigl[x^{\widetilde{M}_j} \bigr]
= \prod_{j=1}^J \frac{1-\tilde{\rho}_j}{1-\tilde{\rho}_j x}
\]
for any $x > 0$ with $\rho x < 1$ {(note that $\tilde{\rho}_j x<1$ for
all $j$ if and only if $\rho x = \max_j \tilde{\rho}_j x<1$, by Lemma
\ref{lemsys-ind-load})}.
Therefore, for all $\theta> 0$ such that $\rho x = \rho(e^{\theta} -
1)/\theta< 1$, we have
\[ \label{eqworkload-tail-upper1}
\limsup_{L\rightarrow\infty} \frac{1}{L} \log\Prob\Biggl(\sum
_{\ell=1}^M U_{\ell} \geq L \Biggr)
\leq\limsup_{L\rightarrow
\infty} \frac{1}{L} \log\Biggl\{ \exp(-\theta
L) \prod_{j=1}^J \frac{1-\tilde{\rho}_j}{1-\tilde{\rho}_j x} \Biggr
\} = -\theta.
\]

Taking the infimum over $\theta$ satisfying $\rho(e^{\theta} -
1)/\theta< 1$, we have established~\eqref{eqtail-workload-upper},
that is,
\[
\limsup_{L \to\infty} \frac{1}{L} \log\Prob\Biggl(\sum
_{\ell
=1}^M U_{\ell} \geq L \Biggr)
\leq-\theta^*.
\]
We now prove the converse inequality.
\iffalse
%
%eB.15 #&#
\begin{equation}
\label{eqtail-workload-lower} \liminf_{L \to\infty} \frac{1}{L} \log
\Prob
\Biggl(\sum_{\ell
=1}^M U_{\ell} \geq
L \Biggr) \geq-\theta^*.
\end{equation}
%
\fi
Without loss of generality, suppose that $\rho= \tilde{\rho}_1$,
and $\widetilde{M}_1$ is a geometric random variable with parameter
$1-\rho$.
Then we can couple $\sum_{\ell=1}^M U_{\ell}$ and $\sum_{\ell
=1}^{\widetilde{M}_1} U_{\ell}$
on the same probability space so that
$\sum_{\ell=1}^M U_{\ell} \geq\sum_{\ell=1}^{\widetilde{M}_1} U_{\ell
}$ with probability $1$.
Thus, it suffices to show that
\[
\liminf_{L \to\infty} \frac{1}{L} \log\Prob\Biggl(\sum
_{\ell
=1}^{\widetilde{M}_1} U_{\ell} \geq L \Biggr)
\geq-\theta^*.
\]
Instead of calculating the quantity directly, consider a $M/D/1$ queue
with load $\rho$,
under the processor-sharing (PS) policy.
Note that for this queuing system, SFA coincides with the PS policy.
By Theorem~\ref{thmSFA1},
$\sum_{\ell=1}^{\widetilde{M}_1} U_{\ell}$ is the steady-state distribution
of the total residual workload in the system.
On the other hand, consider the same queuing system under a FIFO policy.
Since the workload is the same under any work-conserving policy,
$\sum_{\ell=1}^{\widetilde{M}_1} U_{\ell}$ is also the steady-state distribution
of the total workload in this\vspace*{1pt} system, which we denote by $W_{\mathrm{FIFO}}$.
By Theorem 1.4 of~\cite{bigQ}, we can characterize
$\frac{1}{L} \log\Prob(W_{\mathrm{FIFO}} \geq L )$ as follows.\vadjust{\goodbreak}
Let $f(\theta) = \log\bbE[e^{\theta X} ]$,
where $X$ is a Poisson random variable with parameter $\rho$.
Then we have
\[
\lim_{L \to\infty} \frac{1}{L} \log\Prob(W_{\mathrm{FIFO}}
\geq L ) = -\theta^*,
\]
where $\theta^* = \sup\{\theta> 0\dvtx  f(\theta) < \theta\}$.
It is a simple calculation to see that $f(\theta) = \rho(e^{\theta}
- 1)$,
so $\theta^*>0$ satisfies $f(\theta^*) = \theta^*$.
With this lower bound and the upper bound \eqref
{eqtail-workload-upper}, we have established \eqref{eqtail-workload}.
\end{pf*}
%
%{End revised.}

We now provide justifications for Theorem~\ref{thmSFA1}. Consider a
bandwidth-sharing network model as described in Section~\ref{secinsensitive}. Instead of having packets requiring a unit amount of
service, suppose each route $i$ packet has a service requirement that
is independent identically distributed with distribution $\mu_i$ and
mean $1$.
We note that such bandwidth-sharing networks are a special case of the
processor-sharing (PS) queuing network model, as considered by Zachary
\cite{zachary}. In particular, a bandwidth-sharing network is a
processor-sharing network, where network jobs depart the network after
completing service.
General, insensitivity results for the bandwidth-sharing networks
follow as a consequence of the work of Zachary \cite{zachary}.

% \paragraph{Processor-sharing (PS) networks.}
% Consider an open multi-class queuing network.
% Let there be $N$ classes, labeled $1, 2, \ldots, N$.
% For class $i$, \emph{packets} arrive as an independent
% Poisson process of rate $\lambda_i$.
% We also make the following assumption on the service distributions
% requested by incoming packets. Arriving packets of class-$i$
% request services that are i.i.d random variables
% with a common distribution $\mu_i$ and mean $1$, and
% are independent of the service requirements of packets in other
%classes.
% Once a packet receives all its service, it departs the network.
% We denote the number of class-$i$ packets at time $t$ by $M_i(t)$,
% and define the queue-size vector at time $t$ by $\bM(t) = \left(M_i(t)
% If the current queue-size vector is $\mathbf{m}$, then
% class-$i$ packets are assigned a total service rate $\phi_i(
% so that each class-$i$ packet in the system receives a service rate $
% if $m_i > 0$. If $m_i = 0$, then we require that $\phi_i(\mathbf{m})
%= 0$.
% Since the service rate allocation $\bphi(\mathbf{m}) = (\phi_i(
% only depends on the current queue-size vector, as in Section
% the Markov description of the system is given by $\bX(t)$,
% which consists of the queue-size vector $\bM(t) (= \mathbf{m})$ at
%time $t$,
% together with the residual workloads of the set of packets in each
%class.

%On average, $\rho_i = \nu_i/\mu_i$ units of work arrive to class $i$
%per unit time.
%Hence we define $\rho_i = \nu_i/\mu_i$ to be the traffic intensity in
%class $i$.

Following Zachary \cite{zachary}, for $i \in\{1, 2, \ldots, N\}$,
we define the probability distribution $\bar{\mu}_i$
to be the \emph{stationary residual life distribution} of
the renewal process with inter-event distribution $\mu_i$.
That is, if $\mu_i$ has cumulative distribution function $F$,
then $\bar{\mu}_i$ has distribution function $G$ given by
\[
G(x) = 1 - \int_x^{\infty} \bigl(1 - F(y)\bigr)\,dy,\qquad
x\geq0.
\]
Note that if the service requests are deterministically $1$,
that is, $\mu_i$ is the distribution of the deterministic constant $1$,
then $\bar{\mu}_i$ is a uniform distribution on $[0,1]$,
for all $i \in\{1, 2, \ldots, N\}$.

\subsection*{Insensitive rate allocation} Consider a bandwidth-sharing
network described above, with
rate allocation $\bphi(\cdot)$.
A Markovian description of the system is given by a process
$\bY(t)$ which contains the queue-size vector $\bM(t)$ along
with the residual workloads of the set of packets on each route.
If the Markov process $\bY(t)$
admits an invariant measure, then it induces an invariant measure
${\bolds{\pi}}
$ on the process $\bM(t)$.
Such ${\bolds{\pi}}$, when it exists, is called \textit{insensitive} if it depends
on the statistics of the arrivals and service requests only
through the parameters $\bolds{\lambda}= (\lambda_i
)_{i=1}^N$;
in particular, it does not depend
on the detailed service distributions of incoming packets.
A~rate allocation $\bphi(\cdot) = (\phi_i(\cdot) )_{i=1}^N$
is called \textit{insensitive} if it induces an insensitive invariant
measure ${\bolds{\pi}}$ on $\bM(t)$.

It turns out that if the rate allocation $\bphi$ satisfies a \textit{balance property},
then it is insensitive.
%
%deB.1 #&#
\begin{definition}[(Definition 1, \cite{bonaldproutiere2})]
\label{dfbalance}
Consider the bandwidth-sharing network just described.
The rate allocation $\bphi(\cdot)$ is \emph{balanced} if there
exists a function
$\Phi\dvtx  \Integers^{N} \rightarrow\RealsP$ with $\Phi(\mathbf
{0})=1$, and $\Phi(\mathbf{m}) = 0$
for all $\mathbf{m}\notin\IntegersP^{N}$, such that
%
%eB.16 #&#
\begin{equation}
\label{eqbalance} \phi_i(\mathbf{m}) = \frac{\Phi(\mathbf{m}- \be
_i)}{\Phi(\mathbf{m})}\qquad\mbox{for
all } \mathbf{m}\in\IntegersP^{N}, i\in\{1,2,\ldots,N\}.
\end{equation}
\end{definition}

Bonald and Prouti\'{e}re \cite{bonaldproutiere1} proved that
a balanced rate allocation is insensitive with respect to
all phase-type service distributions.
%{DELETE AND MERGE WITH NEXT THEOREM?
%Consider the processor-sharing network described earlier.
%Suppose that all service distributions are phase-type, and
%suppose that the rate allocation $\bphi(\cdot)$ satisfies (
%with $\Phi$ defined in Definition~\ref{dfbalance}.
%Then an invariant measure ${\bolds{\pi}}$ of $\bM(t)$ is given by
%{\bolds{\pi}}(\mathbf{m}) = \Phi(\mathbf{m})\prod_{i=1}^N
%Furthermore, it is unique up to multiplication by a constant factor.
%}
Zachary \cite{zachary} showed that a balanced rate allocation is indeed
insensitive with respect to all general service distributions.
He also gave the characterization of the distribution of the residual workloads
in steady state.

%
%thB.2 #&#
\begin{theorem}[(Theorem 2, \cite{zachary})]\label{thminsengen}
Consider the bandwidth-sharing network described earlier.
A measure ${\bolds{\pi}}$ on $\IntegersP^{N}$ is stationary for $\bM(t)$
and is insensitive to all service distributions with mean $1$,
if and only if it is related to the rate allocation $\bphi$ as follows:
%
%eB.17 #&#
\begin{equation}
\label{eqpbalance} {\bolds{\pi}}(\mathbf{m})\phi_i(\mathbf{m}) = {
\bolds{\pi}}(\mathbf{m}-\be_i)\lambda_i\qquad\mbox{for
all } \mathbf{m}\in\IntegersP^{N},  i \in\{1,2,\ldots,N\},
\end{equation}
where we set ${\bolds{\pi}}(\mathbf{m}-\be_i)$ to be $0$, if $m_i = 0$.
Consequently, ${\bolds{\pi}}$ is given by expression
%
%eB.18 #&#
\begin{equation}
\label{eqmeasure} {\bolds{\pi}}(\mathbf{m}) = \Phi(\mathbf{m})\prod
_{i=1}^N \lambda_i^{m_i}.
\end{equation}
Furthermore, if ${\bolds{\pi}}$ can be normalized to a probability
distribution,
then $\bY(t)$ is positive recurrent, and in steady state,
the residual workload of each route-$i$ packet
in the network is distributed as $\bar{\mu}_i$,
independent from ${\bolds{\pi}}$,
and is conditionally independent from the residual workloads of other
packets, when we condition on the number of packets on each route of
the network.
\end{theorem}

Note that conditions (\ref{eqbalance}) and (\ref{eqpbalance}) are equivalent.
Suppose that $\bphi(\cdot)$ satisfies~(\ref{eqbalance}),
then an invariant measure ${\bolds{\pi}}$ is given by (\ref{eqmeasure}).
Substituting equation (\ref{eqmeasure}) into equation (\ref
{eqbalance}) gives equation (\ref{eqpbalance}).
Conversely, if equation (\ref{eqpbalance}) is satisfied,
then we can just set $\Phi(\mathbf{m}) = {\bolds{\pi}}(\mathbf{m})/\prod_{i=1}^N
\lambda_i^{m_i}$,
and equations (\ref{eqbalance}) and (\ref{eqmeasure}) are satisfied.

\begin{pf*}{Proof of Theorem~\ref{thmSFA1}}
Theorem~\ref{thmSFA1} is now a fairly easy consequence of
Theorem~\ref{thminsengen}.
Consider a bandwidth-sharing network described in Section~\ref{secinsensitive}.
%It can be regarded as an instance of the PS network model described
%earlier.
The additional structures are the additional capacity constraints
\eqref{eqfeasible},
and that arriving packets only require an unit amount of service,
deterministically.
The capacity constraints \eqref{eqfeasible} impose the necessary
condition for stability,
given by \eqref{eqcapfl}. Recall that all arrival rate vectors
$\bolds{\lambda}$ that
satisfy $\bR\bolds{\lambda}< \bC$ are called \emph{strictly
admissible}.

Consider the bandwith vector $\bphi$ as defined by \eqref{dfsfa1} and
\eqref{dfsfa2}.
As remarked earlier, $\bphi$ is admissible, that is, it satisfies the
capacity constraints
\eqref{eqfeasible}. It is balanced by definition, and hence is insensitive
by Theorem~\ref{thminsengen}.
Thus, it induces an stationary measure ${\bolds{\pi}}$ on the
queue-size vector
$\bM(t)$,
given by \eqref{eqmeasure}. For a strictly admissible arrival rate
vector $\bolds{\lambda}$,
the measure is finite, with the normalizing constant $\Phi$ given by
\eqref{eqnormalizeSFA}.
Hence, we can normalize ${\bolds{\pi}}$ to obtain the unique stationary
probability distribution for $\bM(t)$.\vadjust{\goodbreak}

Finally, using Theorem~\ref{thminsengen} and the fact that
all service requests are deterministically $1$,
we see that the stationary residual workloads
are all uniformly distributed on $[0,1]$ and independent.
\end{pf*}\vspace*{-6pt}

%sC #&#
\section{Proof of Lemma \texorpdfstring{\protect\ref{LEMMONOTONE}}{5.1}}

Let $\bD\geq\bZero$, and let
$(\tilde{\alpha}_{\bolds{\sigma}})_{\bolds{\sigma}\in\sS}$ be an
optimal solution to
the program $\PRIMAL(\bD)$.
Then $\tilde{\alpha}_{\bolds{\sigma}} \geq0$ for all $\bolds{\sigma}\in
\sS$,
$\bD\leq\sum_{\bolds{\sigma}\in\sS} \tilde{\alpha}_{\bolds{\sigma}}
\bolds{\sigma}$,
and $\sum_{\bolds{\sigma}\in\sS} \tilde{\alpha}_{\bolds{\sigma}} = \rho
(\bD)$.
We will construct $(\alpha_{\bolds{\sigma}})_{\bolds{\sigma}\in\sS}$
from $(\tilde{\alpha}_{\bolds{\sigma}})_{\bolds{\sigma}\in\sS}$
such that $\alpha_{\bolds{\sigma}} \geq0$ for all $\bolds{\sigma}\in\sS$,
$\bD= \sum_{\bolds{\sigma}\in\sS} \alpha_{\bolds{\sigma}} \bolds{\sigma}$
and $\sum_{\bolds{\sigma}\in\sS} \alpha_{\bolds{\sigma}} = \rho(\bD)$.

If $\bD= \sum_{\bolds{\sigma}\in\sS} \tilde{\alpha}_{\bolds{\sigma}}
\bolds{\sigma}$, then
there is nothing to prove. Thus, suppose that there exists $i$ such that
$D_i < \sum_{\bolds{\sigma}\in\sS} \tilde{\alpha}_{\bolds{\sigma}}
\sigma_i$.
\iffalse
We introduce an optimization problem $\PRIMAL'(\bolds{\lambda})$,
which is similar to $\PRIMAL(\bolds{\lambda})$, and which is
defined to be
%
%eC.3 #&#
%eC.4 #&#
%eC.5 #&#
\begin{eqnarray}
\label{defprimal} &&\mbox{minimize } \sum_{\bolds{\sigma}\in\sS}
\alpha_{\bolds{\sigma}}
\\
&&\mbox{subject to } \bolds{\lambda}= \sum_{\bolds{\sigma}\in\sS}
\alpha_{\bolds{\sigma}} \bolds{\sigma}, \label{constreqprimal}
\\
&&\alpha_{\bolds{\sigma}}\in\RealsP\qquad\mbox{for all }\bolds
{\sigma} \in
\sS.
\end{eqnarray}

Clearly, a solution of the $\PRIMAL'(\bD)$ is a feasible solution for
$\PRIMAL(\bD)$. Therefore, to prove the lemma, it is sufficient
to find $ (\alpha^*_{\bolds{\sigma}} )_{\bolds{\sigma}\in\sS}$ that is
an optimal
solution for $\PRIMAL(\bD)$ and satisfies $\sum_{\bolds{\sigma}\in\sS}
\alpha^*_{\bolds{\sigma}}\bolds{\sigma}= \bD$.

Let $ (\alpha'_{\bolds{\sigma}} )_{\bolds{\sigma}\in\sS}$ be an
optimal solution to
$\PRIMAL(\bD)$. Then
%
\[
\sum_{\bolds{\sigma}\in\sS} \alpha'_{\bolds{\sigma}}
\bolds{\sigma}\geq\bD.
\]
%
If all the inequality constraints are tight, then there is nothing to
prove. Therefore,
suppose that
%
\[
\theta_i \equiv\sum_{\bolds{\sigma}\in\sS}
\alpha'_{\bolds{\sigma}}\sigma_i > D_i
\]
%
for some $i \in\{1, 2, \ldots, N\}$. \fi
We now modify
$ (\tilde{\alpha}_{\bolds{\sigma}} )_{\bolds{\sigma}\in\sS}$ to reduce the
``gap'' between
$\sum_{\bolds{\sigma}\in\sS} \tilde{\alpha}_{\bolds{\sigma}}\sigma_i$
and $D_i$.
%{DELETE??: That is,
%we find some $\left(\alpha''_{\bolds{\sigma}}\right)_{\bolds{\sigma}
%$$\sum_{\bolds{\sigma}\in\sS} \alpha'_{\bolds{\sigma}}\sigma_k =
%for all $k \neq i$, and
%$$\sum_{\bolds{\sigma}\in\sS} \alpha'_{\bolds{\sigma}}\sigma_i >
%and that
%$$\sum_{\bolds{\sigma}\in\sS} \alpha'_{\bolds{\sigma}} = \sum_{
%}

Indeed, since $\sum_{\bolds{\sigma}\in\sS} \tilde{\alpha}_{\bolds{\sigma
}}\sigma_i >
D_i \geq0$,
there is some $\tilde{\bolds{\sigma}} \in\sS$ such that $\tilde{\sigma
}_i = 1$,
and $\tilde{\alpha}_{\tilde{\bolds{\sigma}}} >0$.
We now modify $(\tilde{\alpha}_{\bolds{\sigma}})$ by reducing
$\tilde{\alpha}_{\tilde{\bolds{\sigma}}}$ by a positive amount
\[
\beps= \min(\tilde{\alpha}_{\tilde{\bolds{\sigma}}}, \theta_i -
D_i ),
\]
increasing $\tilde{\alpha}_{\tilde{\bolds{\sigma}} - \be_i}$ by $\beps> 0$,
and keeping all other $\tilde{\alpha}_{\bolds{\sigma}}$ the same
($\tilde{\bolds{\sigma}} - \be_i \in\sS$ by Assumption~\ref{assmonotone}).
Then it is easy to check that $\sum_{\bolds{\sigma}\in\sS} \tilde{\alpha
}_{\bolds{\sigma}}\sigma_i - D_i$
is reduced by $\beps$, $\sum_{\bolds{\sigma}\in\sS} \tilde{\alpha
}_{\bolds{\sigma}
}\sigma_{\ell} - D_{\ell}$
remains the same for all $\ell\neq i$,
$(\tilde{\alpha}_{\bolds{\sigma}})_{\bolds{\sigma}\in\sS}$ remain nonnegative
and we still have $\sum_{\bolds{\sigma}\in\sS} \tilde{\alpha}_{\bolds
{\sigma}} = \rho
(\bD)$.

\iffalse
Now consider the schedule $\tilde{\bolds{\sigma}} - \be_i$,
where $\be_i$ is the $i$th standard unit vector.
By Assumption~\ref{assmonotone}, $\tilde{\bolds{\sigma}} - \be_i \in\sS$.
Let $\beps= \min(\tilde{\alpha}_{\tilde{\bolds{\sigma}}}, \theta_i -
D_i )$ and define
$ (\alpha'_{\bolds{\sigma}} )_{\bolds{\sigma}\in\sS}$ to be
%
\begin{eqnarray*}
 \alpha'_{\bolds{\sigma}} &=& \tilde{\alpha}_{\bolds{\sigma}}\qquad\mbox
{if }\bolds{\sigma} \neq\tilde{\bolds{\sigma}}, \tilde{\bolds{\sigma}} -
\be_i,
\\
 \alpha'_{\tilde{\bolds{\sigma}}} &=& \tilde{\alpha}_{\tilde{\bolds
{\sigma}}} -
\beps\quad\mbox{and}\quad\alpha'_{\tilde{\bolds{\sigma}}-\be_i} =
\tilde{
\alpha}_{\tilde{\bolds{\sigma}}-\be_i} + \beps. %\equiv\sum_{\bolds{
\end{eqnarray*}
%
It then follows that
%
\[
\sum_{\bolds{\sigma}\in\sS} \alpha'_{\bolds{\sigma}}
\sigma_i = \sum_{\bolds{\sigma}\in
\sS} \tilde{
\alpha}_{\bolds{\sigma}} \sigma_i - \beps
\]
%
and
%
\[
\sum_{\bolds{\sigma}\in\sS} \alpha'_{\bolds{\sigma}}
\sigma_{\ell} = \sum_{\bolds{\sigma}
\in\sS} \tilde{
\alpha}_{\bolds{\sigma}} \sigma_{\ell} \mbox{ for all } \ell\neq i.
\]
%
We also have $\alpha'_{\bolds{\sigma}} \geq0$ for all $\bolds{\sigma}\in
\sS$,
$\sum_{\bolds{\sigma}\in\sS} \alpha'_{\bolds{\sigma}}\bolds{\sigma}\geq
\bD$, and $\sum_{\bolds{\sigma}} \alpha'_{\bolds{\sigma}} = \sum_{\bolds
{\sigma}} \alpha
''_{\bolds{\sigma}}$.
\fi
By repeating this procedure
finitely many times, it follows that we can modify $(\tilde{\alpha
}_{\bolds{\sigma}})_{\bolds{\sigma}\in\sS}$
to make it satisfy \eqref{eqcvx}.
This completes the proof of Lemma~\ref{LEMMONOTONE}.

%sD #&#
\section{Proof of Lemma \texorpdfstring{\protect\ref{LEMPOSREC}}{5.10}}\label{apdxposrec}
First we note that by Theorem~\ref{thmSFA1},
$\bBN$ is positive recurrent under the SFA policy, if $\rho
(\bolds{\lambda}) < 1$.
Starting from any initial state,
it also has a strictly positive probability of
reaching the null-state $(\bM(\cdot), \bmu(\cdot)) = \bZero$
at some finite time.
Since the evolution of the virtual system $\bBN$ does not depend on that
of $\bSN$, it is, on its own, positive recurrent.
Next we argue the positive recurrence of the entire network
state building upon this property of $\bBN$.

Sufficient conditions to establish positive Harris recurrence and
ergodicity of a discrete-time
Markov chain $\bX(\tau)$ with state space $\sX$ are given by the
following (see, \cite{Ass}, pages 198--202, and
\cite{Foss-Fluid}, Section~4.2, for details):
\begin{longlist}[(C2)]
\item[(C1)] There exists {a bounded set} $A \in\cB_\sX$ such that
%
%eD.1 #&#
%eD.2 #&#
\begin{eqnarray}
\Prob_{\mathbf{x}} ( T_{A} < \infty) & = & 1 \qquad\mbox{for
any $\mathbf{x}\in\sX$} \label{eqsf1}
\\
\sup_{\mathbf{x}\in A} \bbE_\mathbf{x} [ T_{A} ] & <
& \infty. \label{eqsf2}
\end{eqnarray}
In the above, the stopping time $T_A = \inf\{\tau\geq1 \dvtx  \bX(\tau)
\in A\}$; notation
$\Prob_{\mathbf{x}}(\cdot) \equiv\Prob( \cdot| \bX(0) = \mathbf
{x})$ and
$\bbE_{\mathbf{x}} [\cdot] \equiv\bbE[\cdot| \bX
(0) = \mathbf{x} ]$.

\item[(C2)] Given $A$ satisfying \eqref{eqsf1}--\eqref{eqsf2},
there exists $\mathbf{x}^* \in\sX$, finite $\ell\geq1$ and $\delta
> 0$
such that
%
%eD.3 #&#
%eD.4 #&#
\begin{eqnarray}
{\Prob}_{\mathbf{x}} \bigl(\bX(\ell) = \mathbf{x}^* \bigr) & \geq
&\delta\qquad\mbox{for any $\mathbf{x}\in A$}\label{eqsf3}
\\
{\Prob}_{\mathbf{x}^*} \bigl(\bX(1) = \mathbf{x}^* \bigr) & > & 0.
\label{eqsf4}
\end{eqnarray}
%
%%%\iffalse
%%%\green{
%%%%
%%%\item[{\bf (C2).}] Given $A$ satisfying \eqref{eqsf1}-\eqref{eqsf2},
%%%there exists a probability measure $\lambda(\cdot)$ on $\mathcal
%%%{B}_\sX$, a finite $\ell\geq1$ and $\delta> 0$
%%%such that, for each $\mathbf{x}\in A$,
%%%%
%%%\begin{eqnarray}
%%%{\Prob}_{\mathbf{x}}\big(\bX(\ell) \in B) & \geq&\delta\lambda
%%%(B), \qquad\mbox{for all $B \in\mathcal{B}_\sX$}\label{eqsf3}
%%%\end{eqnarray}
%%%%
%%%\item[{\bf C3.}] The stopping time $T_A$ is aperiodic when $\bX(\cdot
%%%)$ is started from initial distribution $\lambda(\cdot)$.
%%%}
%%%\fi
\end{longlist}
%
%Given this, the following path-wise ergodic property is satisfied (cf.
%any $\bx\in\sX$ and nonnegative measurable function $f\dvtx  \sX\to
%$$ \lim_{T\to\infty} \frac{1}{T} \sum_{\tau= 0}^{T-1} f(X(\tau))
%Here $\E_\pi[f] = \int f(z) \pi(z).$ Note that $\E_\pi[f]$ may
%not be finite.
\old{\green{???????????????????????????????
Conditions \eqref{eqsf1} and \eqref{eqsf3} ensure that $A$ is a
regeneration set and thus $\bX(\cdot)$ is a Harris chain.
Conditions \eqref{eqsf1} and {(C2)} ensure a stationary measure
exists, \cite{Ass}, Theorem 3.2. Condition \eqref{eqsf1}, then,
ensures the stationary measure can be renormalized as a stationary
probability measure. Following this, \eqref{eqsf4} guarantees our
Markov chain is aperiodic thus is ergodic \cite{Ass}, Theorem 3.6. As
we remarked before Lemma~\ref{LEMPOSREC}, the state-space $\sX$ may
have points which are not neighbourhood recurrent. Even so, our
stationary distribution is well defined and we can characterize the set
of recurrent sets in $\sX$ as exactly the sets with positive measure
with respect to our stationary distribution \cite{Ass}, Corollary 3.3.}}

Next, we verify conditions {(C1)} and {(C2)}. For the set of
points where $\bBN$ is empty, say $A$,  condition {(C1)}
follows immediately from the following facts:  (a)~$\bBN$ is
positive recurrent and hence $(\bM(\cdot), \bmu(\cdot))$ returns
to $\bZero$ state in finite expected time starting from any finite
state; (b) $\bD(\cdot)$ is always bounded due to Lemma~\ref{lemdom};
and (c) $\bQ(\cdot)$ returns to the bounded set $\sum_i Q_i(\cdot)
\leq K (N+2)$
whenever $\bM(\cdot) = \bZero$ due to Lemma~\ref{lemdom1}.
Condition {(C2)} can be verified for the null-state $\mathbf{x}^* =
\bZero$
as follows: (a) $(\bM(\cdot), \bmu(\cdot))$ returns to the null
state with positive
probability; (b) given this, it remains there for further $K(N+2) + 1$
time with strictly
positive probability due to Poisson arrival process; (c) in this
additional time $K(N+2) +1$,
the $\bQ(\cdot)$ and $\bD(\cdot)$ are driven to $\bZero$. To see
(c), observe that
when $\bM(\cdot) = \bZero$, $\bD(\cdot) \in\IntegersP^N$. By
construction of our
policy and Assumption~\ref{assmonotone} on structure of $\sS$, it
follows that if
$\bM(\cdot)$ continues to remain $\bZero$, the $\sum_i D_i(\cdot)$
is reduced by at least
unit amount till $\bD(\cdot) = \bZero$; at which moment $\bQ(\cdot
)$ reaches $\bZero$
as well. Since $\sum_i D_i(\cdot) \leq K(N+2)$ by Lemma~\ref{lemdom}, it follows
that $\bM(\cdot)$ need to remain $\bZero$ for this to happen only
for $K(N+2) + 1$
amount of time. This completes the verification of the conditions {(C1)} and {(C2)}.
Consequently, we establish that the network Markov chain, represented
by $\bX(\cdot)$,
is positive recurrent and ergodic.

\end{appendix}

\section*{Acknowledgments} Devavrat Shah and Yuan Zhong would like
to thank John Tsitsiklis for a careful reading of the paper which has
helped improve
the readability, and for his insights and support of this project.
We would also like to thank the referee and Associated Editor
for many valuable comments, which helped with resolving some technical issues
and improving the readability.

%suskaldyti doi

% imsref loaded by linak, 2014-02-20 15:57:55
% imsref loaded by linak, 2014-02-20 16:21:06
%
% imsref loaded by linak, 2014-04-29 13:39:34

\printaddresses


\begin{thebibliography}{41}

%b1 #&#
\bibitem{Ass}
%
\begin{bbook}[mr]
\bauthor{\bsnm{Asmussen},~\bfnm{S{\o}ren}\binits{S.}}
(\byear{2003}).
\btitle{Applied Probability and Queues},
\bedition{2nd} ed.
\bpublisher{Springer},
\blocation{New York}.
\bid{mr={1978607}}
\end{bbook}
%
\bptok{imsref}%
% NOT OUTPUTED:
% isbn = 0-387-00211-1
% fpage = xii+438
\endbibitem

%b2 #&#
\bibitem{BrK}
%
\begin{barticle}[mr]
\bauthor{\bsnm{Birkhoff},~\bfnm{Garrett}\binits{G.}}
(\byear{1946}).
\btitle{Three observations on linear algebra}.
\bjournal{Univ. Nac. Tucum\'an. Revista A.}
\bvolume{5}
\bpages{147--151}.
\bid{mr={0020547}}
\end{barticle}
%
\bptok{imsref}%
\endbibitem

%b3 #&#
\bibitem{bonaldproutiere1}
%
\begin{barticle}[auto:STB|2014/02/12|12:18:25]
\bauthor{\bsnm{Bonald},~\bfnm{T.}\binits{T.}} \AND
\bauthor{\bsnm{Prouti{\`e}re},~\bfnm{A.}\binits{A.}}
(\byear{2002}).
\btitle{Insensitivity in processor-sharing networks}.
\bjournal{Perform. Eval.}
\bvolume{49}
\bpages{193--209}.
\end{barticle}
%
\bptok{imsref}%
% NOT OUTPUTED:
% number = 1-4
\endbibitem

%b4 #&#
\bibitem{bonaldproutiere2}
%
\begin{barticle}[mr]
\bauthor{\bsnm{Bonald},~\bfnm{T.}\binits{T.}} \AND
\bauthor{\bsnm{Prouti{\`e}re},~\bfnm{A.}\binits{A.}}
(\byear{2003}).
\btitle{Insensitive bandwidth sharing in data networks}.
\bjournal{Queueing Syst.}
\bvolume{44}
\bpages{69--100}.
\bid{doi={10.1023/A:1024094807532}, issn={0257-0130}, mr={1989867}}
\end{barticle}
%
\bptok{imsref}%
% NOT OUTPUTED:
% issn = 0257-0130
% url = http://dx.doi.org/10.1023/A:1024094807532
% number = 1
% fjournal = Queueing Systems. Theory and Applications
\endbibitem

%b5 #&#
\bibitem{bramson}
%
\begin{barticle}[mr]
\bauthor{\bsnm{Bramson},~\bfnm{Maury}\binits{M.}}
(\byear{1998}).
\btitle{State space collapse with application to heavy traffic limits
for multiclass queueing networks}.
\bjournal{Queueing Syst.}
\bvolume{30}
\bpages{89--148}.
\bid{doi={10.1023/A:1019160803783}, issn={0257-0130}, mr={1663763}}
\end{barticle}
%
\bptok{imsref}%
% NOT OUTPUTED:
% issn = 0257-0130
% url = http://dx.doi.org/10.1023/A:1019160803783
% number = 1-2
% fjournal = Queueing Systems. Theory and Applications
\endbibitem

%b6 #&#
\bibitem{LD05}
%
\begin{barticle}[mr]
\bauthor{\bsnm{Dai},~\bfnm{J.~G.}\binits{J.~G.}} \AND
\bauthor{\bsnm{Lin},~\bfnm{Wuqin}\binits{W.}}
(\byear{2005}).
\btitle{Maximum pressure policies in stochastic processing networks}.
\bjournal{Oper. Res.}
\bvolume{53}
\bpages{197--218}.
\bid{doi={10.1287/opre.1040.0170}, issn={0030-364X}, mr={2131925}}
\end{barticle}
%
\bptok{imsref}\vadjust{\goodbreak}%
% NOT OUTPUTED:
% issn = 0030-364X
% url = http://dx.doi.org/10.1287/opre.1040.0170
% number = 2
% coden = OPREAI
% fjournal = Operations Research
\endbibitem

%b7 #&#
\bibitem{LD08}
%
\begin{barticle}[mr]
\bauthor{\bsnm{Dai},~\bfnm{J.~G.}\binits{J.~G.}} \AND
\bauthor{\bsnm{Lin},~\bfnm{Wuqin}\binits{W.}}
(\byear{2008}).
\btitle{Asymptotic optimality of maximum pressure policies in
stochastic processing networks}.
\bjournal{Ann. Appl. Probab.}
\bvolume{18}
\bpages{2239--2299}.
\bid{doi={10.1214/08-AAP522}, issn={1050-5164}, mr={2473656}}
\end{barticle}
%
\bptok{imsref}%
% NOT OUTPUTED:
% issn = 1050-5164
% url = http://dx.doi.org/10.1214/08-AAP522
% number = 6
% fjournal = The Annals of Applied Probability
\endbibitem

%b8 #&#
\bibitem{daibala}
%
\begin{bincollection}[auto:STB|2014/02/12|12:18:25]
\bauthor{\bsnm{Dai},~\bfnm{J.~G.}\binits{J.~G.}} \AND
\bauthor{\bsnm{Prabhakar},~\bfnm{B.}\binits{B.}}
(\byear{2000}).
\btitle{The throughput of switches with and without speed-up}.
In \bbooktitle{Proceedings of IEEE Infocom}
\bpages{556--564}.
\bpublisher{IEEE}, \blocation{New York}.
\end{bincollection}
%
\bptok{imsref}%
\endbibitem

%b9 #&#
\bibitem{DKS90}
%
\begin{barticle}[auto:STB|2014/02/12|12:18:25]
\bauthor{\bsnm{Demers},~\bfnm{A.}\binits{A.}},
\bauthor{\bsnm{Keshav},~\bfnm{S.}\binits{S.}} \AND
\bauthor{\bsnm{Shenker},~\bfnm{S.}\binits{S.}}
(\byear{1990}).
\btitle{Analysis and simulation of a fair queuing algorithm}.
\bjournal{Internetworking: Research and Experience}
\bvolume{1}
\bpages{3--26}.
\end{barticle}
%
\bptok{imsref}%
\endbibitem

%b10 #&#
\bibitem{EMPS}
%
\begin{barticle}[mr]
\bauthor{\bparticle{El}~\bsnm{Gamal},~\bfnm{Abbas}\binits{A.}},
\bauthor{\bsnm{Mammen},~\bfnm{James}\binits{J.}},
\bauthor{\bsnm{Prabhakar},~\bfnm{Balaji}\binits{B.}} \AND
\bauthor{\bsnm{Shah},~\bfnm{Devavrat}\binits{D.}}
(\byear{2006}).
\btitle{Optimal throughput-delay scaling in wireless networks. {II}.
{C}onstant-size packets}.
\bjournal{IEEE Trans. Inform. Theory}
\bvolume{52}
\bpages{5111--5116}.
\bid{doi={10.1109/TIT.2006.883548}, issn={0018-9448}, mr={2300378}}
\end{barticle}
%
\bptok{imsref}%
% NOT OUTPUTED:
% issn = 0018-9448
% url = http://dx.doi.org/10.1109/TIT.2006.883548
% number = 11
% coden = IETTAW
% fjournal = Institute of Electrical and Electronics Engineers.
%Transactions on Information Theory
\endbibitem

%b11 #&#
\bibitem{Foss-Fluid}
%
\begin{barticle}[mr]
\bauthor{\bsnm{Foss},~\bfnm{Serguei}\binits{S.}} \AND
\bauthor{\bsnm{Konstantopoulos},~\bfnm{Takis}\binits{T.}}
(\byear{2004}).
\btitle{An overview of some stochastic stability methods}.
\bjournal{J. Oper. Res. Soc. Japan}
\bvolume{47}
\bpages{275--303}.
\bid{issn={0453-4514}, mr={2174067}}
\end{barticle}
%
\bptok{imsref}%
% NOT OUTPUTED:
% issn = 0453-4514
% number = 4
% coden = JORJAS
% fjournal = Journal of the Operations Research Society of Japan
\endbibitem

%b12 #&#
\bibitem{GP93}
%
\begin{barticle}[auto:STB|2014/02/12|12:18:25]
\bauthor{\bsnm{Gallager},~\bfnm{R.}\binits{R.}} \AND
\bauthor{\bsnm{Parekh},~\bfnm{A.}\binits{A.}}
(\byear{1993}).
\btitle{A generalized processor sharing approach to flow control in
integrated services networks: The single-node case}.
\bjournal{IEEE/ACM Transactions on Networking}
\bvolume{1}
\bpages{344--357}.
\end{barticle}
%
\bptok{imsref}%
% NOT OUTPUTED:
% number = 3
\endbibitem

%b13 #&#
\bibitem{bigQ}
%
\begin{bbook}[mr]
\bauthor{\bsnm{Ganesh},~\bfnm{Ayalvadi}\binits{A.}},
\bauthor{\bsnm{O'Connell},~\bfnm{Neil}\binits{N.}} \AND
\bauthor{\bsnm{Wischik},~\bfnm{Damon}\binits{D.}}
(\byear{2004}).
\btitle{Big Queues}.
\bpublisher{Springer},
\blocation{Berlin}.
\bid{doi={10.1007/b95197}, mr={2045489}}
\end{bbook}
%
\bptok{imsref}%
% NOT OUTPUTED:
% isbn = 3-540-20912-3
% url = http://dx.doi.org/10.1007/b95197
% fpage = xii+254
\endbibitem

%b14 #&#
\bibitem{mike2}
%
\begin{bincollection}[mr]
\bauthor{\bsnm{Harrison},~\bfnm{J.~Michael}\binits{J.~M.}}
(\byear{1995}).
\btitle{Balanced fluid models of multiclass queueing networks: A~heavy
traffic conjecture}.
In \bbooktitle{Stochastic Networks}.
\bseries{IMA Vol. Math. Appl.}
\bvolume{71}
\bpages{1--20}.
\bpublisher{Springer},
\blocation{New York}.
\bid{mr={1381003}}
\end{bincollection}
%
\bptok{imsref}%
\endbibitem

%b15 #&#
\bibitem{harrisoncanonical}
%
\begin{barticle}[mr]
\bauthor{\bsnm{Harrison},~\bfnm{J.~Michael}\binits{J.~M.}}
(\byear{2000}).
\btitle{Brownian models of open processing networks: Canonical
representation of workload}.
\bjournal{Ann. Appl. Probab.}
\bvolume{10}
\bpages{75--103}.
\bid{doi={10.1214/aoap/1019737665}, issn={1050-5164}, mr={1765204}}
\end{barticle}
%
\bptok{imsref}%
% NOT OUTPUTED:
% issn = 1050-5164
% url = http://dx.doi.org/10.1214/aoap/1019737665
% number = 1
% fjournal = The Annals of Applied Probability
\endbibitem

%b16 #&#
\bibitem{harrisoncanonicalcorr}
%
\begin{barticle}[mr]
\bauthor{\bsnm{Harrison},~\bfnm{J.~Michael}\binits{J.~M.}}
(\byear{2003}).
\btitle{Correction: ``{B}rownian models of open processing networks:
Canonical representation of workload'' [\textit{Ann. Appl. Probab.}
\textbf{10} (2000), no. 1, 75--103; {MR}1765204 (2001g:60230)]}.
\bjournal{Ann. Appl. Probab.}
\bvolume{13}
\bpages{390--393}.
\bid{doi={10.1214/aoap/1042765673}, issn={1050-5164}, mr={1952004}}
\end{barticle}
%
\bptok{imsref}%
% NOT OUTPUTED:
% issn = 1050-5164
% url = http://dx.doi.org/10.1214/aoap/1042765673
% number = 1
% fjournal = The Annals of Applied Probability
\endbibitem

%b17 #&#
\bibitem{JS}
%
\begin{bincollection}[auto:STB|2014/02/12|12:18:25]
\bauthor{\bsnm{Jagabathula},~\bfnm{S.}\binits{S.}} \AND
\bauthor{\bsnm{Shah},~\bfnm{D.}\binits{D.}}
(\byear{2008}).
\btitle{Optimal delay scheduling in networks with arbitrary constraints}.
In \bbooktitle{Proceedings of the 2008 ACM SIGMETRICS International
Conference on Measurement and Modeling of Computer Systems}
\bpages{395--406}.
\bpublisher{ACM},
\blocation{New York}.\vadjust{\goodbreak}
\end{bincollection}
%
\bptok{imsref}%
\endbibitem

%b18 #&#
\bibitem{kelly-williamsssc}
%
\begin{barticle}[mr]
\bauthor{\bsnm{Kang},~\bfnm{W.~N.}\binits{W.~N.}},
\bauthor{\bsnm{Kelly},~\bfnm{F.~P.}\binits{F.~P.}},
\bauthor{\bsnm{Lee},~\bfnm{N.~H.}\binits{N.~H.}} \AND
\bauthor{\bsnm{Williams},~\bfnm{R.~J.}\binits{R.~J.}}
(\byear{2009}).
\btitle{State space collapse and diffusion approximation for a network
operating under a fair bandwidth sharing policy}.
\bjournal{Ann. Appl. Probab.}
\bvolume{19}
\bpages{1719--1780}.
\bid{doi={10.1214/08-AAP591}, issn={1050-5164}, mr={2569806}}
\end{barticle}
%
\bptok{imsref}%
% NOT OUTPUTED:
% issn = 1050-5164
% url = http://dx.doi.org/10.1214/08-AAP591
% number = 5
% fjournal = The Annals of Applied Probability
\endbibitem

%b19 #&#
\bibitem{Ke79}
%
\begin{bbook}[mr]
\bauthor{\bsnm{Kelly},~\bfnm{Frank~P.}\binits{F.~P.}}
(\byear{1979}).
\btitle{Reversibility and Stochastic Networks}.
\bpublisher{Wiley},
\blocation{Chichester}.
\bid{mr={0554920}}
\end{bbook}
%
\bptok{imsref}%
% NOT OUTPUTED:
% isbn = 0-471-27601-4
% fpage = viii+230
\endbibitem

%b20 #&#
\bibitem{KMW}
%
\begin{barticle}[mr]
\bauthor{\bsnm{Kelly},~\bfnm{F.~P.}\binits{F.~P.}},
\bauthor{\bsnm{Massouli{\'e}},~\bfnm{L.}\binits{L.}} \AND
\bauthor{\bsnm{Walton},~\bfnm{N.~S.}\binits{N.~S.}}
(\byear{2009}).
\btitle{Resource pooling in congested networks: Proportional fairness
and product form}.
\bjournal{Queueing Syst.}
\bvolume{63}
\bpages{165--194}.
\bid{doi={10.1007/s11134-009-9143-8}, issn={0257-0130}, mr={2576010}}
\end{barticle}
%
\bptok{imsref}%
% NOT OUTPUTED:
% issn = 0257-0130
% url = http://dx.doi.org/10.1007/s11134-009-9143-8
% number = 1-4
% fjournal = Queueing Systems. Theory and Applications
\endbibitem

%b21 #&#
\bibitem{kingmanht}
%
\begin{barticle}[mr]
\bauthor{\bsnm{Kingman},~\bfnm{J.~F.~C.}\binits{J.~F.~C.}}
(\byear{1962}).
\btitle{On queues in heavy traffic}.
\bjournal{J. Roy. Statist. Soc. Ser. B}
\bvolume{24}
\bpages{383--392}.
\bid{issn={0035-9246}, mr={0148146}}
\end{barticle}
%
\bptok{imsref}%
% NOT OUTPUTED:
% issn = 0035-9246
% fjournal = Journal of the Royal Statistical Society. Series B.
%Methodological
\endbibitem

%b22 #&#
\bibitem{Matousek2002}
%
\begin{bbook}[mr]
\bauthor{\bsnm{Matou{\v{s}}ek},~\bfnm{Ji{\v{r}}{\'{\i}}}\binits{J.}}
(\byear{2002}).
\btitle{Lectures on Discrete Geometry}.
\bpublisher{Springer},
\blocation{New York}.
\bid{doi={10.1007/978-1-4613-0039-7}, mr={1899299}}
\end{bbook}
%
\bptok{imsref}%
% NOT OUTPUTED:
% isbn = 0-387-95373-6
% url = http://dx.doi.org/10.1007/978-1-4613-0039-7
% fpage = xvi+481
\endbibitem

%b23 #&#
\bibitem{Meyn08}
%
\begin{barticle}[mr]
\bauthor{\bsnm{Meyn},~\bfnm{Sean}\binits{S.}}
(\byear{2009}).
\btitle{Stability and asymptotic optimality of generalized maxweight policies}.
\bjournal{SIAM J. Control Optim.}
\bvolume{47}
\bpages{3259--3294}.
\bid{doi={10.1137/06067746X}, issn={0363-0129}, mr={2476439}}
\bptnote{check year}%
\end{barticle}
%
\bptok{imsref}%
% NOT OUTPUTED:
% issn = 0363-0129
% url = http://dx.doi.org/10.1137/06067746X
% number = 6
% fjournal = SIAM Journal on Control and Optimization
\endbibitem

%b24 #&#
\bibitem{MT93}
%
\begin{bbook}[mr]
\bauthor{\bsnm{Meyn},~\bfnm{S.~P.}\binits{S.~P.}} \AND
\bauthor{\bsnm{Tweedie},~\bfnm{R.~L.}\binits{R.~L.}}
(\byear{1993}).
\btitle{Markov Chains and Stochastic Stability}.
\bpublisher{Springer},
\blocation{London}.
\bid{mr={1287609}}
\end{bbook}
%
\bptok{imsref}%
% NOT OUTPUTED:
% isbn = 3-540-19832-6
% fpage = xvi+ 548
\endbibitem

%b25 #&#
\bibitem{NeelyModiano}
%
\begin{barticle}[auto:STB|2014/02/12|12:18:25]
\bauthor{\bsnm{Neely},~\bfnm{M.}\binits{M.}},
\bauthor{\bsnm{Modiano},~\bfnm{E.}\binits{E.}} \AND
\bauthor{\bsnm{Cheng},~\bfnm{Y.}\binits{Y.}}
(\byear{2007}).
\btitle{Logarithmic delay for $n\times n$ packet switches under the
crossbar constraint}.
\bjournal{IEEE/ACM Transactions on Networking (TON)}
\bvolume{15}
\bpages{657--668}.
\end{barticle}
%
\bptok{imsref}%
% NOT OUTPUTED:
% number = 3
\endbibitem

%b26 #&#
\bibitem{PTh}
%
\begin{bmisc}[auto:STB|2014/02/12|12:18:25]
\bauthor{\bsnm{Prouti{\`e}re},~\bfnm{A.}\binits{A.}}
(\byear{2003}).
\btitle{Insensitivity and stochastic bounds in queuing
netwokrs---applications to flow level traffic modelling in
telecommunication networks.
Ph.D. thesis, Ecole Doctorale, de l'Ecole Polytechnique.}
\end{bmisc}
%
\bptok{imsref}%
\endbibitem

%b27 #&#
\bibitem{STZ25}
%
\begin{bmisc}[auto:STB|2014/02/12|12:18:25]
\bauthor{\bsnm{Shah},~\bfnm{D.}\binits{D.}},
\bauthor{\bsnm{Tsitsiklis},~\bfnm{J.~N.}\binits{J.~N.}} \AND
\bauthor{\bsnm{Zhong},~\bfnm{Y.}\binits{Y.}}
\bhowpublished{Queue-size scaling for input-queued switches.
Unpublished manuscript}.
\end{bmisc}
%
\bptok{imsref}%
% NOT OUTPUTED:
% sortkey = Shah
\endbibitem

%b28 #&#
\bibitem{STZSIGM}
%
\begin{bincollection}[auto:STB|2014/02/12|12:18:25]
\bauthor{\bsnm{Shah},~\bfnm{D.}\binits{D.}},
\bauthor{\bsnm{Tsitsiklis},~\bfnm{J.~N.}\binits{J.~N.}} \AND
\bauthor{\bsnm{Zhong},~\bfnm{Y.}\binits{Y.}}
(\byear{2010}).
\btitle{Qualitative properties of $\alpha$-weighted scheduling policies}.
In \bbooktitle{ACM SIGMETRICS Performance Evaluation Review}
\bpages{239--250}.
\bpublisher{ACM},
\blocation{New York}.
\end{bincollection}
%
\bptok{imsref}%
\endbibitem

%b29 #&#
\bibitem{STZopen}
%
\begin{barticle}[mr]
\bauthor{\bsnm{Shah},~\bfnm{Devavrat}\binits{D.}},
\bauthor{\bsnm{Tsitsiklis},~\bfnm{John~N.}\binits{J.~N.}} \AND
\bauthor{\bsnm{Zhong},~\bfnm{Yuan}\binits{Y.}}
(\byear{2011}).
\btitle{Optimal scaling of average queue sizes in an input-queued
switch: An open problem}.
\bjournal{Queueing Syst.}
\bvolume{68}
\bpages{375--384}.
\bid{doi={10.1007/s11134-011-9234-1}, issn={0257-0130}, mr={2834209}}
\end{barticle}
%
\bptok{imsref}%
% NOT OUTPUTED:
% issn = 0257-0130
% url = http://dx.doi.org/10.1007/s11134-011-9234-1
% number = 3-4
% fjournal = Queueing Systems. Theory and Applications
\endbibitem

%b30 #&#
\bibitem{STZ}
%
\begin{barticle}[mr]
\bauthor{\bsnm{Shah},~\bfnm{D.}\binits{D.}},
\bauthor{\bsnm{Tsitsiklis},~\bfnm{J.~N.}\binits{J.~N.}} \AND
\bauthor{\bsnm{Zhong},~\bfnm{Y.}\binits{Y.}}
(\byear{2014}).
\btitle{Qualitative properties of {$\alpha$}-fair policies in
bandwidth-sharing networks}.
\bjournal{Ann. Appl. Probab.}
\bvolume{24}
\bpages{76--113}.
\bid{doi={10.1214/12-AAP915}, issn={1050-5164}, mr={3161642}}
\bptnote{check year}%
\end{barticle}
%
\bptok{imsref}%
% NOT OUTPUTED:
% issn = 1050-5164
% url = http://dx.doi.org/10.1214/12-AAP915
% number = 1
% fjournal = The Annals of Applied Probability
\endbibitem

%b31 #&#
\bibitem{SWo}
%
\begin{barticle}[mr]
\bauthor{\bsnm{Shah},~\bfnm{Devavrat}\binits{D.}} \AND
\bauthor{\bsnm{Wischik},~\bfnm{Damon}\binits{D.}}
(\byear{2011}).
\btitle{Fluid models of congestion collapse in overloaded switched networks}.
\bjournal{Queueing Syst.}
\bvolume{69}
\bpages{121--143}.
\bid{doi={10.1007/s11134-011-9250-1}, issn={0257-0130}, mr={2836736}}
\end{barticle}
%
\bptok{imsref}%
% NOT OUTPUTED:
% issn = 0257-0130
% url = http://dx.doi.org/10.1007/s11134-011-9250-1
% number = 2
% fjournal = Queueing Systems. Theory and Applications
\endbibitem

%b32 #&#
\bibitem{SW}
%
\begin{barticle}[mr]
\bauthor{\bsnm{Shah},~\bfnm{Devavrat}\binits{D.}} \AND
\bauthor{\bsnm{Wischik},~\bfnm{Damon}\binits{D.}}
(\byear{2012}).
\btitle{Switched networks with maximum weight policies: Fluid
approximation and multiplicative state space collapse}.
\bjournal{Ann. Appl. Probab.}
\bvolume{22}
\bpages{70--127}.
\bid{doi={10.1214/11-AAP759}, issn={1050-5164}, mr={2932543}}
\end{barticle}
%
\bptok{imsref}%
% NOT OUTPUTED:
% issn = 1050-5164
% url = http://dx.doi.org/10.1214/11-AAP759
% number = 1
% fjournal = The Annals of Applied Probability
\endbibitem

%b33 #&#
\bibitem{stolyar}
%
\begin{barticle}[mr]
\bauthor{\bsnm{Stolyar},~\bfnm{Alexander~L.}\binits{A.~L.}}
(\byear{2004}).
\btitle{Maxweight scheduling in a generalized switch: State space
collapse and workload minimization in heavy traffic}.
\bjournal{Ann. Appl. Probab.}
\bvolume{14}
\bpages{1--53}.
\bid{doi={10.1214/aoap/1075828046}, issn={1050-5164}, mr={2023015}}
\end{barticle}
%
\bptok{imsref}%
% NOT OUTPUTED:
% issn = 1050-5164
% url = http://dx.doi.org/10.1214/aoap/1075828046
% number = 1
% fjournal = The Annals of Applied Probability
\endbibitem

%b34 #&#
\bibitem{Stolyar-LDP}
%
\begin{barticle}[mr]
\bauthor{\bsnm{Stolyar},~\bfnm{Alexander~L.}\binits{A.~L.}}
(\byear{2008}).
\btitle{Large deviations of queues sharing a randomly time-varying server}.
\bjournal{Queueing Syst.}
\bvolume{59}
\bpages{1--35}.
\bid{doi={10.1007/s11134-008-9072-y}, issn={0257-0130}, mr={2429895}}
\end{barticle}
%
\bptok{imsref}%
% NOT OUTPUTED:
% issn = 0257-0130
% url = http://dx.doi.org/10.1007/s11134-008-9072-y
% number = 1
% fjournal = Queueing Systems. Theory and Applications
\endbibitem

%b35 #&#
\bibitem{tassiula1}
%
\begin{barticle}[mr]
\bauthor{\bsnm{Tassiulas},~\bfnm{Leandros}\binits{L.}} \AND
\bauthor{\bsnm{Ephremides},~\bfnm{Anthony}\binits{A.}}
(\byear{1992}).
\btitle{Stability properties of constrained queueing systems and
scheduling policies for maximum throughput in multihop radio networks}.
\bjournal{IEEE Trans. Automat. Control}
\bvolume{37}
\bpages{1936--1948}.
\bid{doi={10.1109/9.182479}, issn={0018-9286}, mr={1200609}}
\end{barticle}
%
\bptok{imsref}%
% NOT OUTPUTED:
% issn = 0018-9286
% url = http://dx.doi.org/10.1109/9.182479
% number = 12
% coden = IETAA9
% fjournal = Institute of Electrical and Electronics Engineers.
%Transactions on Automatic Control
\endbibitem

%b36 #&#
\bibitem{VL-LDP}
%
\begin{bincollection}[auto:STB|2014/02/12|12:18:25]
\bauthor{\bsnm{Venkataramanan},~\bfnm{V.~J.}\binits{V.~J.}} \AND
\bauthor{\bsnm{Lin},~\bfnm{X.}\binits{X.}}
(\byear{2007}).
\btitle{Structural properties of LDP for queue-length based wireless
scheduling algorithms}.
In \bbooktitle{45th Annual Allerton Conference on Communication, Control, and Computing}
\bpages{759--766}.
\bpublisher{IEEE}, \blocation{New York}.
\end{bincollection}
%
\bptok{imsref}%
\endbibitem

%b37 #&#
\bibitem{VN}
%
\begin{bincollection}[mr]
\bauthor{\bparticle{von} \bsnm{Neumann},~\bfnm{John}\binits{J.}}
(\byear{1953}).
\btitle{A certain zero-sum two-person game equivalent to the optimal
assignment problem}.
In \bbooktitle{Contributions to the Theory of Games}.
\bpages{5--12}.
\bpublisher{Princeton Univ. Press},
\blocation{Princeton, NJ}.
\bid{mr={0054920}}
\end{bincollection}
%
\bptok{imsref}%
\endbibitem

%b38 #&#
\bibitem{walton}
%
\begin{barticle}[mr]
\bauthor{\bsnm{Walton},~\bfnm{N.~S.}\binits{N.~S.}}
(\byear{2009}).
\btitle{Proportional fairness and its relationship with multi-class
queueing networks}.
\bjournal{Ann. Appl. Probab.}
\bvolume{19}
\bpages{2301--2333}.
\bid{doi={10.1214/09-AAP612}, issn={1050-5164}, mr={2588246}}
\end{barticle}
%
\bptok{imsref}%
% NOT OUTPUTED:
% issn = 1050-5164
% url = http://dx.doi.org/10.1214/09-AAP612
% number = 6
% fjournal = The Annals of Applied Probability
\endbibitem

%b39 #&#
\bibitem{whittspl}
%
\begin{bbook}[mr]
\bauthor{\bsnm{Whitt},~\bfnm{Ward}\binits{W.}}
(\byear{2002}).
\btitle{Stochastic-Process Limits}.
\bpublisher{Springer},
\blocation{New York}.
%application to queues}.
\bid{mr={1876437}}
\bptnote{check year}%
\end{bbook}
%
\bptok{imsref}%
% NOT OUTPUTED:
% isbn = 0-387-95358-2
% fpage = xxiv+602
\endbibitem

%b40 #&#
\bibitem{williams}
%
\begin{barticle}[mr]
\bauthor{\bsnm{Williams},~\bfnm{R.~J.}\binits{R.~J.}}
(\byear{1998}).
\btitle{Diffusion approximations for open multiclass queueing
networks: Sufficient conditions involving state space collapse}.
\bjournal{Queueing Syst.}
\bvolume{30}
\bpages{27--88}.
\bid{doi={10.1023/A:1019108819713}, issn={0257-0130}, mr={1663759}}
\end{barticle}
%
\bptok{imsref}%
% NOT OUTPUTED:
% issn = 0257-0130
% url = http://dx.doi.org/10.1023/A:1019108819713
% number = 1-2
% fjournal = Queueing Systems. Theory and Applications
\endbibitem

%b41 #&#
\bibitem{zachary}
%
\begin{barticle}[mr]
\bauthor{\bsnm{Zachary},~\bfnm{Stan}\binits{S.}}
(\byear{2007}).
\btitle{A note on insensitivity in stochastic networks}.
\bjournal{J. Appl. Probab.}
\bvolume{44}
\bpages{238--248}.
\bid{doi={10.1239/jap/1175267175}, issn={0021-9002}, mr={2312999}}
\end{barticle}
%
\bptok{imsref}%
% NOT OUTPUTED:
% issn = 0021-9002
% url = http://dx.doi.org/10.1239/jap/1175267175
% number = 1
% coden = JPRBAM
% fjournal = Journal of Applied Probability
\endbibitem

\end{thebibliography}
\end{document}